\documentclass{amsart}

\usepackage{hyperref}
\usepackage{amsmath}
\usepackage{amssymb}
\usepackage{amscd}
\usepackage{graphicx}
\usepackage{mathdots}
\usepackage{gensymb}
\usepackage{epstopdf}
 \usepackage{todonotes}
 \usepackage{youngtab}
 \usepackage{wrapfig}

\usepackage{epsfig}       
\usepackage{epic,eepic}        

\newtheorem{thm}{Theorem}[section]
\newtheorem{prop}[thm]{Proposition}
\newtheorem{lem}[thm]{Lemma}
\newtheorem{cor}[thm]{Corollary}
\newtheorem{conj}[thm]{Conjecture}

\theoremstyle{definition}

\theoremstyle{remark}

\numberwithin{equation}{section}

\newcommand{\group}[1]{\mathrm{#1}}

\newcommand{\mon}{\vec{H}}

\newcommand{\Tr}{\operatorname{Tr}}

\newcommand{\D}{\mathbb{D}}

\renewcommand{\P}{\mathbb{P}}

\newcommand{\LIS}{\mathrm{LIS}}

\newcommand{\C}{\mathbb{C}}

\begin{document}


\title{On the Complex Asymptotics of the HCIZ and BGW Integrals}

 
\author{Jonathan Novak}
\address{Department of Mathematics, University of California, San Diego, USA}
\email{jinovak@ucsd.edu}



 \maketitle
 
 \begin{quote}
 	\textit{There seems to be a connection between large $N$ and the permutation groups.
	--- Stuart Samuel, 1980}
\end{quote}
 
 \tableofcontents
 \pagebreak




\section{Introduction}

\subsection{Objective}
The purpose of this paper is to prove a longstanding conjecture on the $N \to \infty$ asymptotic behavior
of the Harish-Chandra/Itzykson-Zuber (HCIZ) integral,

	\begin{equation*}
		\label{eqn:HCIZ} 
		I_N = \int_{\mathrm{U}(N)} e^{zN \mathrm{Tr} AUBU^{-1}} \mathrm{d}U,
	\end{equation*}
	
\noindent
and its additive counterpart, the Br\'ezin-Gross-Witten (BGW) integral,

	\begin{equation*}
		\label{eqn:BGW}
		J_N = \int_{\mathrm{U}(N)} e^{zN \mathrm{Tr}(AU + BU^{-1})} \mathrm{d}U.
	\end{equation*}
	
\noindent
These are integrals over $N \times N$ unitary matrices against unit mass Haar measure,
the integrands of which depend on a complex parameter $z$ and a pair of $N \times N$ complex matrices
$A$ and $B$. The conjecture we prove emerged from a cluster of 1980 theoretical physics papers on the 
large $N$ limit of $\group{U}(N)$ lattice gauge theory \cite{Bars,BG,GW,IZ,Samuel,Wadia}, 
and has been of perennial interest in physics ever since; see the reviews \cite{BBMP,Morozov,ZZ}. 
It entered mathematics in the early 2000s along with growing interest in random matrices, and 
was precisely formulated in work of Collins \cite[Section 5]{Collins}, Guionnet \cite[Section 4.3]{Guionnet:PS}, and
Zelditch \cite[Section 4]{Z}. The conjecture has since attained the status of an outstanding open problem in asymptotic analysis, 
and has become perhaps the most prominent question at the confluence of random matrix theory and representation 
theory; see e.g. \cite{BGH} for a recent perspective. It may be stated as follows.

Given a Young diagram $\alpha$ with $d$ cells, $\ell(\alpha)$ rows, and 
$\alpha_i$ cells in the $i$th row, let 

	\begin{equation*}
		p_\alpha(x_1,\dots,x_N) = \prod_{i=1}^{\ell(\alpha)} \sum_{j=1}^N x_j^{\alpha_i}
	\end{equation*}
	
\noindent
be the corresponding Newton power sum symmetric polynomial in $N$ variables. 

\begin{conj}
	\label{conj:Main}
	Given any $M \geq 0,$ there exists a corresponding $\varepsilon_M>0$ such that, for any integer $k \geq 0$,
		
		\begin{equation*}
			I_N = e^{\sum_{g=0}^k N^{2-2g}F_N^{(g)}+ o(N^{2-2k})} \quad\text{ and }\quad 
			J_N = e^{\sum_{g=0}^k N^{2-2g}G_N^{(g)} + o(N^{2-2k})}
		\end{equation*}
		
	\noindent
	as $N \to \infty$, where the error term is uniform over complex numbers $z$ of modulus at most 
	$\varepsilon_M$ and complex matrices $A,B$ of spectral radius at most $M$, and
		
	\begin{equation*}
	\begin{split}
			F_N^{(g)} &= \sum_{d=1}^\infty  \frac{z^d}{d!} \sum_{\alpha,\beta \vdash d}
			\frac{p_\alpha(a_1,\dots,a_N)}{N^{\ell(\alpha)}} \frac{p_\beta(b_1,\dots,b_N)}{N^{\ell(\beta)}} F_g(\alpha,\beta), \\
			G_N^{(g)} &= \sum_{d=1}^\infty \frac{z^{2d}}{d!} \sum_{\beta \vdash d}
			\frac{p_\beta(c_1,\dots,c_N)}{N^{\ell(\beta)}} G_g(\beta) ,
	\end{split}
	\end{equation*}
	
	\noindent
	are analytic functions of $z$, the eigenvalues $a_1,\dots,a_N$ of $A$, the eigenvalues $b_1,\dots,b_N$ of $B$,
	and the eigenvalues $c_1,\dots,c_N$ of $C=AB$.
	Moreover, the coefficients $F_g(\alpha,\beta)$ and $G_g(\beta)$ are integers.
\end{conj}

The main result of this paper is a proof of Conjecture \ref{conj:Main}. Before outlining our argument,
let us briefly unpack the conjecture's meaning. Its salient feature is the claim that $I_N$ and $J_N$ admit what 
physicists call ``strong coupling expansions'' --- their logarithms have complete $N \to \infty$ asymptotic expansions on the scale $N^{2-2g}$
provided the ``coupling constant'' $z$ is sufficiently small and the ``external fields'' $A$ and $B$ are uniformly bounded
(our parameter $z$ is inversely proportional to the physical coupling constant, so that small $|z|$ corresponds to strong coupling).
Without loss in generality, we may take $M=1$ as the uniform bound
on the spectral radii of $A$ and $B$. The conjecture then asserts the existence of $\varepsilon>0$ such that, for any
given $k \geq 0,\kappa > 0$, there is a corresponding $N(k,\kappa)$ with $N \geq N(k,\kappa)$ 
implying

	\begin{equation*}
		\left| \log I_N - \sum_{g=0}^k N^{2-2g} F_N^{(g)} \right| \leq \kappa N^{2-2k} \quad\text{ and }\quad
		\left| \log J_N - \sum_{g=0}^k N^{2-2g} G_N^{(g)} \right| \leq \kappa N^{2-2k}
	\end{equation*}
	
\noindent
for all complex numbers $z$ of modulus at most $\varepsilon$ and all complex matrices $A,B$ with eigenvalues of modulus 
at most $1,$ where ``$\log$'' denotes the principal branch of the complex logarithm. The coefficients of these purported
asymptotic expansions --- the  ``free energies'' $F_N^{(g)}$ and $G_N^{(g)}$ --- 
are themselves dependent on $N,$ and hence could conceivably interact with the asymptotic scale. 
The conjecture addresses this by further claiming that $F_N^{(g)}$ and $G_N^{(g)}$ are analytically determined by the data
$(z,A,B)$ in a manner which precludes this possibility: it implies the bounds

	\begin{equation*}
		|F_N^{(g)}| \leq \sum_{d=1}^\infty \frac{\varepsilon^d}{d!} \sum_{\alpha,\beta \vdash d} |F_g(\alpha,\beta)|
		\quad\text{ and }\quad
		|G_N^{(g)}| \leq \sum_{d=1}^\infty \frac{\varepsilon^{2d}}{d!} \sum_{\beta \vdash d} |G_g(\beta)|,
	\end{equation*}
	
\noindent
which are finite and depend only on $\varepsilon$ and $g.$
Finally, the conjecture asserts that the universal coefficients $F_g(\alpha,\beta)$ and $G_g(\beta),$
which determine $F_N^{(g)}$ and $G_N^{(g)}$ but do not depend on the data $(z,A,B),$
are integers. This claim is rooted in the notion of ``topological expansion,''
a fundamental but analytically non-rigorous principle in quantum field theory
which generalizes the apparatus of Feynman diagrams to 
matrix integrals \cite{tHooft1,BIZ,BIPZ,DGZ,IZ,Witten}, and beyond \cite{EO,KS}. This principle predicts that
the structure constants $F_g(\alpha,\beta)$ and $G_g(\beta)$ are combinatorial invariants of 
compact connected genus $g$ Riemann surfaces.

\subsection{Results}
The main result of this paper is a proof of Conjecture \ref{conj:Main}. Our argument proeeds in
three stages: exact formulas, stable asymptotics, and functional asymptotics. 

\subsubsection{Exact formulas}
Our point of departure is a pair of novel absolutely convergent series expansions of
$I_N$ and $J_N$ which are amenable to large $N$ analysis. 

	\begin{thm}
		\label{thm:StringExpansions}
		For any $N \in \mathbb{N}$, we have
		
			\begin{align*}
				I_N &= 1 + \sum_{d=1}^\infty \frac{z^d}{d!} \P(\LIS_d \leq N) \sum_{\alpha,\beta \vdash d} 
				 p_\alpha(a_1,\dots,a_N) p_\beta(b_1,\dots,b_N) 
				 \langle \omega_\alpha \Omega_N^{-1} \omega_\beta \rangle,\\
				J_N &= 1 + \sum_{d=1}^\infty \frac{z^{2d}}{d!} N^d\P(\LIS_d \leq N) \sum_{\beta \vdash d} 
				p_\beta(c_1,\dots,c_N) \langle \Omega_N^{-1} \omega_\beta  \rangle,
			\end{align*}
			
		\noindent
		where $\P(\LIS_d \leq N)$ is the probability that a uniformly random permutation from the symmetric group 
		$\group{S}(d)$ has no increasing subsequence of length $N+1$, and
		$\langle \omega_\alpha \Omega_N^{-1} \omega_\beta \rangle$ is the Plancherel expectation of a certain
		natural observable of Young diagrams with $d$ cells and at most $N$ rows. 
		These series converge absolutely and uniformly on compact subsets of $\C^{2N+1}$ and $\C^{N+1}$, respectively.
	\end{thm}

\noindent
We call these series the ``string expansions'' of $I_N$ and $J_N;$ this terminology is explained in 
Section \ref{sec:Exact} below. In the absence of external fields, the string expansion of the BGW integral reduces to the beautiful formula

	\begin{equation*}
		\int_{\mathrm{U}(N)} e^{zN \mathrm{Tr}(U + U^{-1})} 
		\mathrm{d}U = 1+ \sum_{d=1}^\infty \frac{z^{2d}}{d!} N^{2d}\P(\LIS_d \leq N),
	\end{equation*}
	
\noindent
which is due to Rains \cite{Rains} and equivalent to a result of Gessel \cite{Gessel}.
The Gessel-Rains identity was the starting point of  Baik, Deift, and Johansson \cite{BDJ} 
in their seminal work showing that the $d \to \infty$ fluctuations of $\LIS_d$ around its asymptotic mean value 
of $2\sqrt{d}$ are governed by the Tracy-Widom distribution. Informative expositions of this
landmark result may be found in \cite{AD,Romik,Stanley:ICM}. 
The existence of a connection between the HCIZ integral and increasing subsequences appears to 
have been previously unknown. Since the Fourier transform of any unitarily invariant random matrix 
is a mixture of HCIZ integrals, the HCIZ-LIS connection exposes a new and very direct link between
random matrices and random permutations.

\subsubsection{Stable asymptotics}
In Section \ref{sec:Stable}, we analyze the $N \to \infty$ asymptotics of each fixed string coefficient
of $I_N$ and $J_N,$ i.e. the large $N$ asymptotics of the Plancherel expectation $\langle \omega_\alpha \Omega_N^{-1}
\omega_\beta \rangle$ with fixed $\alpha,\beta \vdash d.$ We show that
$\langle \omega_\alpha \Omega_N^{-1} \omega_\beta \rangle$ admits a convergent asymptotic expansion
on the scale $1/N,$ and that this expansion is a generating function for ``monotone'' walks on the Cayley graph of the 
symmetric group $\group{S}(d)$ with boundary conditions $\alpha,\beta.$ 
Monotone walks are self-interacting trajectories: the future of a monotone walk depends on its past. 
It is a fundamental fact, discovered in \cite{Novak:Banach} and further developed in \cite{MN}, 
that these trajectories play the role of Feynman diagrams 
for integration against Haar measure on the unitary group.

For any fixed $N \in \mathbb{N},$ one can replace the first $N$ string coefficients of 
$I_N$ and $J_N$ with their $1/N$ expansions, but not so for higher terms. 
The issue is conceptually similar to that faced when studying the homotopy 
groups of $\group{U}(N)$, which behave regularly at first but eventually become wild. Topologists 
see past this by studying the stable unitary group $\group{U}$, an $N=\infty$ version of
$\group{U}(N)$ which does not suffer from this defect \cite{Bott}. The price paid is
that $\group{U}$ is not a Lie group, but an infinite-dimensional manifold which is not locally compact. 
In Section \ref{sec:Stable}, we introduce the stable HCIZ and BGW integrals, $I$ and $J$, 
which are $N=\infty$ versions of $I_N$ and $J_N$. Conceptually, these objects are the integrals

	\begin{equation*}
		I = \int_{\group{U}} e^{\frac{z}{\hbar} \Tr AUBU^{-1}} \mathrm{d}U \quad\text{ and }\quad 
		J = \int_{\group{U}} e^{\frac{z}{\hbar} \Tr (AU+BU^{-1})} \mathrm{d}U,
	\end{equation*}
	
\noindent
with $\hbar$ an infinitely small parameter, $A$ and $B$ infinitely large matrices, and $\mathrm{d}U$ the non-existent Haar measure
on the stable unitary group $\group{U}.$  Like the homotopy groups of $\group{U},$
the topological expansions of $I$ and $J$ can be completely understood; the price paid is that $I$ and $J$ are 
not analytic functions, but formal power series in infinitely many variables which are not convergent. 

\begin{thm}
	\label{thm:MainStable}
		We have 
		
			\begin{equation*}
				I = e^{\sum_{g=0}^\infty \hbar^{2g-2} F^{(g)}} \quad\text{ and} \quad
				J = e^{\sum_{g=0}^\infty \hbar^{2g-2}G^{(g)}},
			\end{equation*}
			
		\noindent
		where the stable free energies $F^{(g)}$ and $G^{(g)}$ 
		are generating functions for the genus $g$ monotone double and single Hurwitz numbers, 
		respectively.
\end{thm}
	
Hurwitz theory, familiar to algebraic geometers as the prototypical enumerative theory of maps from curves to 
curves, plays a prominent role in contemporary enumerative geometry; see \cite{ELSV,GJV,KL,OP}, and 
\cite{DYZ} for a recent overview.
Monotone Hurwitz theory \cite{GGN1,GGN2,GGN3,GGN4,GGN5} 
is a desymmetrized version of classical Hurwitz theory which, 
rather surprisingly, is exactly solvable to exactly the same extent. Just as there are explicit formulas for 
classical Hurwitz numbers in genus zero and one \cite{Hurwitz,Vakil}, there are
explicit formulas for monotone Hurwitz numbers in genus zero and one \cite{GGN1,GGN2}, and the two sets of
formulas are structurally analogous. Monotone Hurwitz numbers manifest versions of 
polynomiality \cite{GGN2} and integrability \cite{GGN4} which mirror the polynomiality \cite{ELSV} and integrability \cite{Okounkov:MRL}
of their classical counterparts. The consonance between the classical and monotone theories is to some extent
explained by the fact that both are governed by the Eynard-Orantin topological recursion formalism ---
the two theories are structurally identical, but are generated by different spectral curves \cite{BEMS,DDM}.

Monotone Hurwitz theory has proved to be a useful tool 
with diverse applications \cite{BK,CDO,GGR,Montanaro,Novak:JSP}, and its discovery has sparked a surge of interest
in combinatorial deformations of classical Hurwitz numbers \cite{ACEH,ALS,CD,DDM,DK,DKPS,HKL}.
Although the subject has taken on a life of its own, monotone Hurwtz 
numbers were originally summoned from the void as a weapon with which to 
attack Conjecture \ref{conj:Main}. In this paper, they fulfill their initial purpose.

\subsubsection{Functional asymptotics}
In order to prove Conjecture \ref{conj:Main}, we must descend from the stable world of $N =\infty$ the unstable
world of finite $N.$ To navigate this passage, we must address the questions of convergence and approximation which 
are the analytic substance of Conjecture \ref{conj:Main}. Prior knowledge of the stable limit, which comprises the 
combinatorial substance of Conjecture \ref{conj:Main}, is extremely useful in this regard --- 
since we know what the answer is supposed to be, the analysis becomes a task of verification rather than discovery. 

More precisely, if Conjecture \ref{conj:Main} is true then the free energies $F_N^{(g)}$ and $G_N^{(g)}$
must be generating functions for monotone Hurwitz numbers of genus $g.$ Remarkably, monotone Hurwitz theory guarantees that the stable
free energies $F^{(g)}$ and $G^{(g)}$ remain stable at finite $N$: replacing the formal parameter $\hbar$ with $N^{-1}$
and the formal alphabets $A,B,C$ with the spectra of uniformly bounded $N \times N$ complex matrices yields absolutely summable power 
series. Even better, the radius of convergence of these series is bounded below by a positive constant $\delta$ 
independent of both $N$ and $g.$ We thus have explicit analytic candidates for $F_N^{(g)}$ and $G_N^{(g)},$ 
with a stable domain of holomorphy. 

The stable analyticity of $F_N^{(g)}$ and $G_N^{(g)}$ does not mean that one can deduce Conjecture \ref{conj:Main} from 
Theorem \ref{thm:MainStable} simply by replacing $\hbar$ with $N^{-1}$ --- this fails because the series

	\begin{equation*}
		F_N = \sum_{g=0}^\infty N^{2-2g}F_N^{(g)} \quad\text{ and }\quad G_N = \sum_{g=0}^\infty N^{2-2g}G_N^{(g)}
	\end{equation*}
	
\noindent
are not uniformly convergent on any nondegenerate polydisc for any finite $N.$ 
This is typical of generating functions associated
with 2D quantum gravity \cite{DGZ,Witten}, and one sees similar phenomena in the world of maps on surfaces and 
Hermitian matrix integrals \cite{EM,Maurel}. The divergence of these series forces the introduction of a cutoff at
fixed genus $g=k$, and an ensuing analysis of the holomorphic discrepancy functions

	\begin{equation*}
		1- \frac{I_N}{e^{\sum_{g=0}^k N^{2-2g}F_N^{(g)}}} \quad\text{ and }\quad 
		1-\frac{J_N}{e^{\sum_{g=0}^k N^{2-2g}G_N^{(g)}}}.
	\end{equation*}
	
\noindent
It is here that knowledge of the full string expansions of $I_N$ and $J_N$ at finite $N$ is essential: it leads to a
 ``topological bound'' which controls the moduli of the discrepancy functions on small polydiscs by a 
 quantity of order $N^{2-2k}.$
Complex analytic tools may then be utilized to convert the topological bound into a topological
approximation, replacing a uniform $O$-term with a uniform $o$-term at the logarithmic scale.
The upshot of this analysis is our main theorem, which proves Conjecture \ref{conj:Main}.

	\begin{thm}
	\label{thm:Main}
		Conjecture \ref{conj:Main} is true, and the structure constants $F_g(\alpha,\beta)$ and $G_g(\beta)$
		are given by
		
			\begin{equation*}
				F_g(\alpha,\beta) = (-1)^{\ell(\alpha)+\ell(\beta)} \mon_g(\alpha,\beta) \quad\text{ and }\quad
				G_g(\beta) = (-1)^{d+\ell(\beta)} \mon_g(\beta),
			\end{equation*}
			
		\noindent
		where $\mon_g(\alpha,\beta)$ and $\mon_g(\beta)$ are the monotone double and single 
		Hurwitz numbers of genus $g$.
	\end{thm}

\subsection{Context}
Conjecture \ref{conj:Main} is the subject of a large literature, and many powerful and impressive 
results have previously been obtained. For the HCIZ integral, the main highlight is Guionnet and Zeitouni's large deviation theory proof \cite{GZ} 
of Matytsin's heuristics \cite{Matytsin}, which characterize the leading asymptotics of $I_N$ in terms of the flow of a compressible fluid.
For the BGW integral, one has Johansson's Toeplitz determinat proof \cite{Johansson:MRL} of Gross and Witten's explicit formula \cite{GW} for the 
leading asymptotics of $J_N$ in the absence of external fields, a result which set the stage for the breakthrough work \cite{BDJ}. 
Another powerful technique is the use of Schwinger-Dyson 
``loop'' equations \cite{Guionnet:book} to obtain both the leading \cite{CGM} and sub-leading \cite{GN} 
asymptotics of a large class of unitary matrix
integrals containing the HCIZ and BGW integrals as prototypes. 

The common limitation of these prior works is that they are restricted to real asymptotics: they are obtained under the additional hypothesis that
both the coupling constant and the eigenvalues of the external fields are real. This assumption is required in order to force
the integrands of $I_N$ and $J_N$ to be positive functions on $\group{U}(N),$ so that probabilistic methods
can be applied. Indeed, all previous approaches to Conjecture \ref{conj:Main} are, ultimately,
elaborations of the classical Laplace method for the asymptotic evaluation of real integrals depending on a large real parameter. 
As soon as complex parameters are allowed, $I_N$ and $J_N$ become oscillatory integrals.
The failure of previous works to treat the complex asymptotics of $I_N$ and $J_N$ is not 
just a technical limitation: many if not most situations in which one would like to invoke the conclusion of
Conjecture \ref{conj:Main} involve complex parameters in an essential way. For example, in order to 
analyze the spectral asymptotics of random matrices using characteristic functions, one needs the asymptotics of the 
orbital integral $I_N$ with complex coupling $z=i,$ which were previously inaccessible, except in certain degenerate scaling
limits \cite{GM,OlshVersh}. This is the sole reason that Fourier analysis has not been a viable technique in the asymptotic spectral analsysi
of random matrices. For exactly the same reason, it has not been possible to make direct use of the Harish-Chandra/Kirillov
formula \cite{HC,K} in asymptotic representation theory. The results of this paper clear the way for
a direct and unified approach to asymptotic random matrix theory and asymptotic representation theory based on Fourier analysis.
Our results can moreover be applied to analyze certain asymptotic problems of physical problems interest which
have been mired in confusion for some time \cite{McNov}.

We have taken a conceptual as opposed to computational approach to the asymptotics of
the HCIZ and BGW integrals by first constructing and understanding their $N=\infty$ stable limits 
and then using this insight to build $N \to \infty$ approximations.
A first pass at this was made in \cite{GGN3}, where Goulden, Guay-Paquet and the author
succeeded in obtaining complete asymptotics for each fixed HCIZ string coefficient,
but failed to understand the full string series at finite $N$ and its remarkable connection with longest increasing subsequences
and Plancherel measure. Consequently, \cite{GGN3} failed to bridge 
the infinitely large gap between $N=\infty$ and $N \to \infty.$ Moreover, the fundamental fact that the relationship between the
HCIZ and BGW integrals is precisely the relationship between double and single Hurwitz numbers was not
perceived in \cite{GGN3}, where the BGW integral was not considered. Indeed, prior to the present work, no matrix model for monotone single Hurwitz numbers was
known, and it was an open question to find one \cite{ALS,DDM}. Given that the known matrix model \cite{BEMS} for classical single Hurwitz numbers
is somewhat contrived, it is remarkable that its counterpart for monotone single Hurwitz numbers is given by none other than the BGW integral, 
the basic special function of lattice gauge theory.


\section{Exact Formulas}
\label{sec:Exact}
	
	In this section, we prove Theorem \ref{thm:StringExpansions}, which is the starting point of our analysis.
												
	\subsection{Symmetric polynomials}
	Given a Young diagram $\alpha$, the associated Newton power sum symmetric polynomial $p_\alpha$ in commuting variables 
	$x_1,\dots,x_N$ is 
	
		\begin{equation*}
			p_\alpha(x_1,\dots,x_N) = \prod_{i=1}^{\ell(\alpha)} 
			\sum_{j=1}^N x_j^{\alpha_i}.
		\end{equation*}	
		
	\noindent
	It is a classical result of Newton (see \cite{Macdonald,Stanley:EC2}) that the polynomials
	
		\begin{equation*}
			p_\alpha(x_1,\dots,x_N), \quad \alpha \vdash d,
		\end{equation*}
		
	\noindent
	span the space $\Lambda_N^{(d)}$ of homogeneous degree $d$ symmetric polynomials in $x_1,\dots,x_N$.
	The Newton polynomials interface naturally with analysis: if $a_1,\dots,a_N$
	is a point configuration in $\mathbb{C}$, then normalized power sums evaluated on these points
	are products of moments of the corresponding empirical probability measure $\mu.$ That is, we have
	
		\begin{equation*}
			\frac{p_\alpha(a_1,\dots,a_N)}{N^{\ell(\alpha)}} = 
			\prod_{i=1}^{\ell(\alpha)} \int_{\mathbb{C}} \zeta^{\alpha_i} \mu(\mathrm{d}\zeta).
		\end{equation*}
		
	\noindent
	In particular, the normalized power sums	
	which appear in Conjecture \ref{conj:Main} are products of moments of the empirical eigenvalue distributions
	of the matrices $A,B$, and $C$. The power sums are the preferred basis for coupling expansions
	in lattice gauge theory, where they are referred to as ``string states'' \cite{BT}.
	
	There is another family of symmetric polynomials which play a role in what follows: the Schur polynomials. 
	Given a Young diagram $\lambda$ with $d$ cells, let $(\mathsf{V}^\lambda,R^\lambda)$
	denote the corresponding irreducible complex representation of the symmetric group
	$\group{S}(d)$, and set
	
		\begin{equation*}
			\chi_\alpha(\lambda) = \Tr R^\lambda(\pi),
		\end{equation*}
		
	\noindent
	where $\pi \in \group{S}(d)$ belongs to the conjugacy class $C_\alpha$ of
	permutations of cycle type $\alpha$. The Schur polynomials,
	
		\begin{equation*}
			s_\lambda(x_1,\dots,x_N) = \frac{1}{d!} \sum_{\alpha \vdash d} 
			|C_\alpha| \chi_\alpha(\lambda) p_\alpha(a_1,\dots,a_N),
			\quad \lambda \vdash d, \ell(\lambda) \leq N.
		\end{equation*}
		
	\noindent
	form a basis of $\Lambda_N^{(d)}$, and the expansion of a given 
	Newton polynomial in the  Schur basis is
	
		\begin{equation*}
			p_\alpha(x_1,\dots,x_N) = \sum_{\substack{\lambda \vdash d \\
			\ell(\lambda) \leq N}} \chi_\alpha(\lambda) s_\lambda(x_1,\dots,x_N).
		\end{equation*}
	
	\noindent
	Evaluations of Schur polynomials at complex points also have representation-theoretic meaning:
	they are  irreducible characters of the general linear group $\group{GL}_N(\C).$ More precisely, 
	given a Young diagram $\lambda$ with at most $N$ rows, 
	let $(\mathsf{W}^\lambda,S^\lambda)$ denote the corresponding irreducible polynomial representation of 
	$\group{GL}_N(\mathbb{C})$.
	Then
	
		\begin{equation*}
			\Tr S^\lambda(A) = s_\lambda(a_1,\dots,a_N)
		\end{equation*}
		
	\noindent
	for any $A \in \mathrm{GL}_N(\mathbb{C})$ with eigenvalues $a_1,\dots,a_N$.

	\subsection{Basic integrals}
	Given a symmetric polynomial $f$ in $N$ variables and an $N \times N$ matrix matrix $A$, 
	write $f(A)$ for the evaluation of $f$ on the spectrum of $A$. We shall need the following basic integration formulas, 
	which are well-known manifestations of Schur orthogonality, see e.g. \cite{Macdonald}. We provide a proof for the 
	sake of completeness.
							
		\begin{lem}
			\label{lem:BasicIntegrals}
			For any Young diagrams $\lambda,\mu$ and matrices $A,B \in \mathrm{Mat}_N(\mathbb{C})$, we
			have
			
				\begin{equation*}
					\int_{\group{U}(N)} s_\lambda(AUBU^{-1}) \mathrm{d}U
					= \frac{s_\lambda(A) s_\lambda(B)}{\dim \mathsf{W}^\lambda}
				\end{equation*}
				
			\noindent
			and
			
				\begin{equation*}
					\int_{\group{U}(N)} s_\lambda(AU) s_\mu(BU^{-1}) \mathrm{d}U
					= \delta_{\lambda\mu}\frac{s_\lambda(AB)}{\dim \mathsf{W}^\lambda}.
				\end{equation*}
		\end{lem}
		
		\begin{proof}
			Suppose first that $A,B \in \group{GL}_N(\mathbb{C})$. Then $AUBU^{-1} \in \group{GL}_N(\C)$, 
			and we have
			
				\begin{align*}
					s_\lambda(AUBU^{-1}) &= \Tr S^\lambda(AUBU^{-1}) \\ 
					&= \Tr S^\lambda(A)S^\lambda(U)S^\lambda(B)S^\lambda(U^{-1}) \\
					&= \sum_{i,j,k,l=1}^N 
					S^\lambda(A)_{ij} S^\lambda(U)_{jk} S^\lambda(B)_{kl} S^\lambda(U^{-1})_{li}.
				\end{align*}
				
			\noindent
			Thus 
			
				\begin{equation*}
					\int_{\group{U}(N)} s_\lambda(AUBU^{-1}) \mathrm{d}U
					= \sum_{i,j,k,l=1}^N S^\lambda(A)_{ij} S^\lambda(B)_{kl}
					\int_{\group{U}(N)} S^\lambda(U)_{jk}S^\lambda(U^{-1})_{li} \mathrm{d}U.
				\end{equation*}
				
			\noindent
			By Schur orthogonality for the matrix elements of an irreducible representation, we have
			
				\begin{equation*}
					\int_{\group{U}(N)} S^\lambda(U)_{jk}S^\lambda(U^{-1})_{li} \mathrm{d}U =
					\frac{\delta_{ij}\delta_{kl}}{\dim \mathsf{W}^\lambda},
				\end{equation*}
				
			\noindent
			and hence
			
				\begin{align*}
					\int_{\group{U}(N)} s_\lambda(AUBU^{-1}) \mathrm{d}U
					&= \frac{1}{\dim \mathsf{W}^\lambda} \sum_{i=1}^N S^\lambda(A)_{ii}
					\sum_{k=1}^N S^\lambda(B)_{kk} \\
					&=  \frac{\Tr S^\lambda(A) \Tr S^\lambda(B)}{\dim \mathsf{W}^\lambda} \\
					&= \frac{s_\lambda(A) s_\lambda(B)}{\dim \mathsf{W}^\lambda}.
				\end{align*}
		
			\noindent
			Similarly, if $A,B \in \group{GL}_N(\C)$, then $AU,BU^{-1} \in \group{GL}_N(\C),$ and
			we have
			
				\begin{align*}
					s_\lambda(AU) = \Tr S^\lambda(A)S^\lambda(U) &= \sum_{i,j=1}^N S^\lambda(A)_{ij} S^\lambda(U)_{ji} \\
					s_\mu(BU^{-1}) = \Tr S^\mu(B)S^\mu(U^{-1})
					&=\sum_{k,l=1}^N S^\mu(B)_{kl} S^\mu(U^{-1})_{lk}.
				\end{align*}
				
			\noindent
			Thus 
			
				\begin{equation*}
					\int_{\group{U}(N)} s_\lambda(AU) s_\mu(BU^{-1}) \mathrm{d}U
					= \sum_{i,j,k,l=1}^N S^\lambda(A)_{ij} S^\mu(B)_{kl}
					\int_{\group{U}(N)} S^\lambda(U)_{ji} S^\mu(U^{-1})_{lk} \mathrm{d}U.
				\end{equation*}
				
			\noindent
			By Schur orthogonality for the matrix elements of different irreducible representations,
			
				\begin{equation*}
					\int_{\group{U}(N)} S^\lambda(U)_{ji} S^\mu(U^{-1})_{lk} \mathrm{d}U 
					= \frac{ \delta_{\lambda\mu} \delta_{il}\delta_{jk}}{\dim \mathsf{W}^\lambda},
				\end{equation*}	
				
			\noindent
			and we conclude that
			
				\begin{equation*}
					\int_{\group{U}(N)} s_\lambda(AU) s_\mu(BU^{-1}) \mathrm{d}U
					= \frac{\delta_{\lambda\mu}}{\dim \mathsf{W}^\lambda}  \sum_{i,j=1}^N S^\lambda(A)_{ij} S^\mu(B)_{ji}
					= \delta_{\lambda\mu} \frac{s_\lambda(AB)}{\dim \mathsf{W}^\lambda}.
				\end{equation*}
				
			That these integral evaluations remain valid for arbitrary complex matrices $A$ and $B$ can be seen by
			taking limits. Let $(A_n)_{n=1}^\infty$ and $(B_n)_{n=1}^\infty$ be sequences in $\group{GL}_N(\C)$ such that
							
				\begin{equation*}			
					\lim_{n \to \infty} A_n = A \quad\text{ and }\quad \lim_{n \to \infty} B_n=B,
				\end{equation*}
			
			\noindent	
			and apply the Dominated Convergence Theorem to obtain
			
				\begin{align*}
					\int_{\group{U}(N)} s_\lambda(AUBU^{-1}) \mathrm{d}U 
					&= \int_{\group{U}(N)} \lim_{n \to \infty} s_\lambda(A_nUB_nU^{-1})  \mathrm{d}U \\
					&= \lim_{n \to \infty} \int_{\group{U}(N)} s_\lambda(A_nUB_nU^{-1})  \mathrm{d}U \\
					&= \lim_{n \to \infty} \frac{s_\lambda(A_n) s_\lambda(B_n)}{\dim \mathsf{W}^\lambda} \\
					&= \frac{s_\lambda(A) s_\lambda(B)}{\dim \mathsf{W}^\lambda}
				\end{align*}
				
			\noindent
			and 
			
				\begin{align*}
					\int_{\group{U}(N)} s_\lambda(AU) s_\mu(BU^{-1}) \mathrm{d}U 
					&= \int_{\group{U}(N)} \lim_{n \to \infty} s_\lambda(A_nU) s_\mu(B_nU^{-1})  \mathrm{d}U \\
					&= \lim_{n \to \infty} \int_{\group{U}(N)} s_\lambda(A_nU) s_\mu(B_nU^{-1})  \mathrm{d}U \\
					&= \lim_{n \to \infty} \delta_{\lambda\mu} \frac{s_\lambda(A_nB_n)}{\dim \mathsf{W}^\lambda} \\
					&= \delta_{\lambda\mu} \frac{s_\lambda(AB)}{\dim \mathsf{W}^\lambda}.
				\end{align*}
				
		\end{proof}
					
		\subsection{Character expansions}		
		Lemma \eqref{lem:BasicIntegrals} leads to the following series representations of 
		$I_N$ and $J_N$ in terms of Schur polynomials. Expansions of this sort appear in various forms in the physics literature, 
		and were perhaps first utilized in work of James \cite{James} in multivariate statistics, where $I_N$ and $J_N$ are treated
		as hypergeometric functions with matrix arguments.
				
		\begin{thm}
			\label{thm:CharacterExpansion}
			For any $z \in \mathbb{C}$, and any $A,B \in \mathrm{Mat}_N(\mathbb{C})$, 
			we have
			
				\begin{align*}
				\int_{\group{U}(N)} e^{zN\Tr AUBU^{-1}} \mathrm{d}U &= 1 + \sum_{d=1}^\infty  \frac{z^d}{d!} N^d
				 \sum_{\substack{\lambda \vdash d \\ \ell(\lambda) \leq N}} 
				 s_\lambda(a_1,\dots,a_N)s_\lambda(b_1,\dots,b_N) \frac{\dim \mathsf{V}^\lambda}{\dim \mathsf{W}^\lambda} \\
				\int_{\group{U}(N)} e^{zN\Tr(AU+BU^{-1})} \mathrm{d}U &= 1 + \sum_{d=1}^\infty  \frac{z^{2d}}{d!d!} N^{2d} 
				 \sum_{\substack{\lambda \vdash d \\ \ell(\lambda) \leq N}} 
				s_\lambda(c_1,\dots,c_N) \frac{(\dim \mathsf{V}^\lambda)^2}{\dim \mathsf{W}^\lambda} ,
				\end{align*}
				
			\noindent
			where $a_1,\dots,a_N$ are the eigenvalues of $A$, $b_1,\dots,b_N$ are the eigenvalues of $B$,
			and $c_1,\dots,c_N$ are the eigenvalues of $C=AB$. These series converge absolutely and uniformly
			on compact subsets of $\C^{2N+1}$ and $\C^{N+1},$ respectively.
		\end{thm}
				
		\begin{proof}
		Consider first the HCIZ integral. Differentiating under the integral sign,
		the Maclaurin series of $I_N$ as an entire function of $z$ is
		
			\begin{align*}
				I_N &= \int_{\group{U}(N)} e^{zN p_1(AUBU^{-1})} \mathrm{d}U \\
				&= 1+ \sum_{d=1}^\infty  \frac{z^d}{d!} N^d
				\int_{\group{U}(N)} p_{1^d}(AUBU^{-1}) \mathrm{d}U   \\
				&= 1+\sum_{d=1}^\infty \frac{z^d}{d!}  N^d
				\sum_{\substack{\lambda \vdash d \\ \ell(\lambda) \leq N}}
				(\dim \mathsf{V}^\lambda) \int_{\group{U}(N)} s_\lambda(AUBU^{-1}) \mathrm{d}U \\
				&=  1 + \sum_{d=1}^\infty \frac{z^d}{d!}  N^d
				\sum_{\substack{\lambda \vdash d \\ \ell(\lambda) \leq N}}
				s_\lambda(A) s_\lambda(B)\frac{\dim \mathsf{V}^\lambda}{\dim \mathsf{W}^\lambda},
			\end{align*}
			
		\noindent
		by Lemma \ref{lem:BasicIntegrals}. 
				
		For the BGW integral, we have 
		
			\begin{align*}
				J_N &= \int_{\group{U}(N)} e^{zNp_1(AU)} e^{zNp_1(BU^{-1})} \mathrm{d}U \\
				&=  1 + \sum_{d=1}^\infty \frac{z^{2d}}{d!d!} N^{2d}
				\sum_{\substack{\lambda \vdash d \\ \ell(\lambda) \leq N}} 
				\sum_{\substack{\mu \vdash d \\ \ell(\mu) \leq N}} 
				(\dim \mathsf{V}^\lambda)(\dim V^\mu)
				\int_{\group{U}(N)} s_\lambda(AU) s_\mu(BU^{-1}) \mathrm{d}U   \\
				&= 1 + \sum_{d=1}^\infty \frac{z^{2d}}{d!d!} N^{2d}
				\sum_{\substack{\lambda \vdash d \\ \ell(\lambda) \leq N}} 
				s_\lambda(AB)\frac{(\dim \mathsf{V}^\lambda)^2}{\dim \mathsf{W}^\lambda},
			\end{align*}
			
		\noindent
		by Lemma \ref{lem:BasicIntegrals}. 		
	\end{proof}
	
	Let us perform a consistency check by examining these formulas in the absence of external fields,
	i.e. when both $A$ and $B$ are the identity matrix.
	For the HCIZ integral, we see directly from the definition that $I_N=e^{zN^2}$ when 
	$A$ and $B$ are the identity. In this case the character expansion of $I_N$ becomes  
	
		\begin{equation*}
			I_N = 1 + \sum_{d=1}^\infty \frac{z^d}{d!}  N^d
				\sum_{\substack{\lambda \vdash d \\ \ell(\lambda) \leq N}}
				(\dim \mathsf{V}^\lambda) (\dim \mathsf{W}^\lambda) = 
				1 + \sum_{d=1}^\infty \frac{z^d}{d!}  N^{2d},
		\end{equation*}
		
	\noindent
	where we have used the isotypic decomposition of the space of $N$-dimensional tensors of rank $d$
	as an $\group{S}(d) \times \group{GL}_N(\C)$ module,
	
		\begin{equation*}
			\left( \mathbb{C}^N \right)^{\otimes d} \simeq \bigoplus_{\substack{\lambda \vdash d \\
			\ell(\lambda) \leq N}} \mathsf{V}^\lambda \otimes \mathsf{W}^\lambda.
		\end{equation*}
		
	\noindent
	For the BGW integral, in the case $AB=I$ the character expansion becomes 
	
		\begin{equation*}
			J_N = 1 + \sum_{d=1}^\infty \frac{z^{2d}}{d!d!} N^{2d}
				\sum_{\substack{\lambda \vdash d \\ \ell(\lambda) \leq N}} 
				(\dim \mathsf{V}^\lambda)^2 = 1 + \sum_{d=1}^\infty \frac{z^{2d}}{d!} N^{2d}
				\P(\LIS_d \leq N),
		\end{equation*}
		
	\noindent
	where in the final equality we used the Robinshon-Schensted correspondence (see below).
	This is exactly the Gessel-Rains identity. Generalizations of the Gessel-Rains identity to integrals
	over truncated unitary matrices were obtained in \cite{Novak:EJC,Novak:IMRN}, and analogues 
	for the other classical groups may be found in \cite{BR,Rains}.
	
	\subsection{String expansions}
	In order to address Conjecture \ref{conj:Main}, we want expansions of 
	$I_N$ and $J_N$ in terms of Newton polynomials rather than Schur polynomials --- 
	string expansions rather than character expansions. We will now obtain the string expansions
	of $I_N$ and $J_N$ from their character expansions.
	
	As is well-known \cite{Macdonald,Stanley:EC2}, 
	the dimension of $\mathsf{V}^\lambda$ is equal to the number of standard Young
	tableaux of shape $\lambda$. Thus, by the Robinson-Schensted correspondence \cite{Stanley:EC2}, we have
		
		\begin{equation*}
			\sum_{\substack{\lambda \vdash d \\ \ell(\lambda) \leq N}}  (\dim \mathsf{V}^\lambda)^2 = |\group{S}_N(d)|,
		\end{equation*}
		
	\noindent
	where $\group{S}_N(d) \subseteq \group{S}(d)$ is the set of permutations with no increasing subsequence
	of length $N+1$. It follows that
	
		\begin{equation*}
			\lambda \mapsto \frac{(\dim \mathsf{V}^\lambda)^2}{|\group{S}_N(d)|}
		\end{equation*}
		
	\noindent
	is the mass function of a probability measure on the set of Young diagrams with $d$ cells an at most $N$ rows.
	This probability measure is known as the (row-restricted) Plancherel measure, see \cite{Kerov,Romik}.
	We denote expectation with respect to Plancherel measure by angled brackets:
	
		\begin{equation*}
			\langle f \rangle = \sum_{\substack{\lambda \vdash d \\ \ell(\lambda) \leq N}}
			f(\lambda) \frac{(\dim \mathsf{V}^\lambda)^2}{|\group{S}_N(d)|}.
		\end{equation*}

	\noindent
	Note that if $N \geq d$ then the restriction on number of rows is vacuous, and the Plancherel measure is a 
	probability measure on the full set of Young diagrams with $d$ cells whose normalization constant is $|\mathrm{S}(d)|=d!.$
					
	For a Young diagram $\alpha \vdash d$, let us identify the conjugacy class 
	$C_\alpha \subset \group{S}(d)$ with the formal sum of its elements, so that
	$C_\alpha$ becomes a central element in the group algebra $\mathbb{C}\group{S}(d)$.
	By Schur's Lemma, $C_\alpha$ acts as a scalar operator in any irreducible 
	representation $(\mathsf{V}^\lambda,R^\lambda)$ of $\mathbb{C}\group{S}(d)$,
	i.e.
	
		\begin{equation*}
			R^\lambda(C_\alpha) = \omega_\alpha(\lambda) I_{\mathsf{V}^\lambda}
		\end{equation*}
		
	\noindent
	where 
	
		\begin{equation*}
			\omega_\alpha(\lambda) = \frac{|C_\alpha| \chi_\alpha(\lambda)}{\dim \mathsf{V}^\lambda}
		\end{equation*}
		
	\noindent
	and $I_{\mathsf{V}^\lambda} \in \mathrm{End} \mathsf{V}^\lambda$ is the identity operator.	
	
	Let us introduce the positive function $\Omega_N$ on Young diagrams with $d$ cells and at most
	$N$ rows defined by
	
		\begin{equation*}
			\Omega_N(\lambda) = \frac{d!}{N^d} \frac{\dim \mathsf{W}^\lambda}{\dim \mathsf{V}^\lambda}
			= \prod_{\Box \in \lambda} \left( 1 + \frac{c(\Box)}{N} \right).
		\end{equation*}
				
	\noindent
	Here we have used the dimension formulas \cite{Macdonald,Stanley:EC2}
	
		\begin{equation*}
			\dim \mathsf{V}^\lambda = \frac{d!}{\prod_{\Box \in \lambda} h(\Box)}
			\quad\text{ and }\quad
			\dim \mathsf{W}^\lambda = \prod_{\Box \in \lambda} \frac{N + c(\Box)}{h(\Box)},
		\end{equation*}
		
	\noindent
	where  $h(\Box)$ is the hook length of a given cell $\Box \in \lambda$ (number of cells to the right of $\Box$
	plus number of cells below $\Box$ plus one) and $c(\Box)$ is its content (column index less row index), 
	to render $\Omega_N(\lambda)$ as an explicit product. Note that
	
		\begin{equation*}
			\Omega_N^{-1}(\lambda) = \prod_{\Box \in \lambda} \frac{1}{1+\frac{c(\Box)}{N}}
		\end{equation*}
		
	\noindent
	is well-defined and positive since $\ell(\lambda) \leq N.$
	The functions $\Omega_N^{\pm 1}$ seem to be closely related to the 
	``$\Omega$-points'' considered by physicists in the context of gauge/string dualities \cite{BT,CMR,GT1}, but which
	seem not to have been fully understood in that context. In terms of $\Omega_N$, 
	the Schur function expansions of $I_N$ and $J_N$ are
	
		\begin{align*}
			I_N &= 1 + \sum_{d=1}^\infty z^d \sum_{\substack{\lambda \vdash d \\ \ell(\lambda) \leq N}} 
			s_\lambda(A) s_\lambda(B) \Omega_N^{-1}(\lambda) \\
			J_N &=  1 + \sum_{d=1}^\infty \frac{z^{2d}}{d!} N^d \sum_{\substack{\lambda \vdash d \\ \ell(\lambda) \leq N}}
			s_\lambda(C) \Omega_N^{-1}(\lambda)
			\dim \mathsf{V}^\lambda.
		\end{align*}
		
	We now prove Theorem \ref{thm:StringExpansions}, which we restate here using the notation just established.
								
	\begin{thm}
		For any $z \in \mathbb{C}$ and any $A,B \in \mathrm{Mat}_N(\mathbb{C})$, 
		we have
			
				\begin{align*}
				\int_{\group{U}(N)} e^{zN\Tr AUBU^{-1}} \mathrm{d}U &= 1 + \sum_{d=1}^\infty \frac{z^d}{d!} \P(\LIS_d \leq N)
				 \sum_{\alpha,\beta \vdash d} p_\alpha(A)p_\beta(B)
				\langle \omega_\alpha \Omega_N^{-1} \omega_\beta \rangle,\\
				\int_{\group{U}(N)} e^{zN\Tr(AU+BU^{-1})} \mathrm{d}U 
				&= 1 + \sum_{d=1}^\infty  \frac{z^{2d}}{d!} N^d\P(\LIS_d \leq N)
				\sum_{\beta \vdash d} p_\beta(C)\langle \Omega_N^{-1} \omega_\beta \rangle,
				\end{align*}
				
			\noindent
			where $C=AB$. 
	\end{thm}
	
	\begin{proof}	
	For the HCIZ integral, we have
	
			\begin{align*}
				I_N&= 1 + \sum_{d=1}^\infty z^d
				\sum_{\substack{\lambda \vdash d \\ \ell(\lambda) \leq N}}
				s_\lambda(A) s_\lambda(B) \Omega_N^{-1}(\lambda)  \\
				&= 1+ \sum_{d=1}^\infty z^d
				\sum_{\substack{\lambda \vdash d \\ \ell(\lambda) \leq N}}
				\left( \sum_{\alpha \vdash d} \frac{|C_\alpha| \chi_\alpha(\lambda)}{d!} p_\alpha(A) \right)
				\left( \sum_{\beta \vdash d} \frac{|C_\beta| \chi_\beta(\lambda)}{d!} p_\beta(B) \right)  \Omega_N^{-1}(\lambda) \\
				&= 1 + \sum_{d=1}^\infty \frac{z^d}{d!} \frac{|S_N(d)|}{d!} \sum_{\alpha,\beta \vdash d} p_\alpha(A) p_\beta(B)
				\sum_{\substack{\lambda \vdash d \\ \ell(\lambda) \leq N}} 
				\frac{|C_\alpha| \chi_\alpha(\lambda)}{\dim \mathsf{V}^\lambda}\Omega_N^{-1}(\lambda)
				\frac{|C_\beta| \chi_\beta(\lambda)}{\dim \mathsf{V}^\lambda} 
				\frac{(\dim \mathsf{V}^\lambda)^2}{|S_N(d)|}  \\
				&= 1+ \sum_{d=1}^\infty \frac{z^d}{d!}\P(\LIS_d \leq N)\sum_{\alpha,\beta \vdash d} p_\alpha(A) p_\beta(B)
				\left\langle \omega_\alpha \Omega_N^{-1} \omega_\beta \right\rangle.
			\end{align*}
			
		For the BGW integral, we have

			\begin{align*}
				J_N &= 1 + \sum_{d=1}^\infty \frac{z^{2d}}{d!} N^d
				\sum_{\substack{\lambda \vdash d \\ \ell(\lambda) \leq N}} 
				s_\lambda(AB) \Omega_N^{-1}(\lambda) \dim \mathsf{V}^\lambda \\
				&= 1 + \sum_{d=1}^\infty  \frac{z^{2d}}{d!} N^d
				\sum_{\substack{\lambda \vdash d \\ \ell(\lambda) \leq N}} \left(
				\sum_{\beta \vdash d} \frac{|C_\beta|\chi_\beta(\lambda)}{d!}
				p_\beta(AB) \right)
				\Omega_N^{-1}(\lambda) \dim \mathsf{V}^\lambda \\
				&= 1 + \sum_{d=1}^\infty  \frac{z^{2d}}{d!}  N^d \frac{|S_N(d)|}{d!} \sum_{\beta \vdash d}
				\sum_{\substack{\lambda \vdash d \\ \ell(\lambda) \leq N}}
				\Omega_N^{-1}(\lambda)\frac{|C_\beta|\chi_\beta(\lambda)}{\dim \mathsf{V}^\lambda} 
				 \frac{(\dim \mathsf{V}^\lambda)^2}{|S_N(d)|} p_\beta(AB) \\
				&= 1 + \sum_{d=1}^\infty \frac{z^{2d}}{d!} N^d \P(\LIS_d \leq N)  
				\sum_{\beta \vdash d} p_\beta(AB)  \langle \Omega_N^{-1} \omega_\beta \rangle.
			\end{align*}
	\end{proof}
	
	\subsection{Basic bounds}
	Let us write the string expansions of $I_N$ and $J_N$ in the form
	
		\begin{equation*}
			I_N = 1 + \sum_{d=1}^\infty \frac{z^d}{d!}I_N(d) \quad\text{ and }\quad
			J_N = 1 + \sum_{d=1}^\infty \frac{z^{2d}}{d!}J_N(d).
		\end{equation*}
		
	\noindent
	Thus $I_N(d)$ and $J_N(d)$ are the symmetric polynomials 
	
		\begin{align*}
			I_N(d) &= N^d
				 \sum_{\substack{\lambda \vdash d \\ \ell(\lambda) \leq N}} s_\lambda(a_1,\dots,a_N) s_\lambda(b_1,\dots,b_N) 
				 \frac{\dim \mathsf{V}^\lambda}{\dim \mathsf{W}^\lambda} \\
				&= \P(\LIS_d \leq N)\sum_{\alpha,\beta \vdash d} p_\alpha(a_1,\dots,a_N)p_\beta(b_1,\dots,b_N) \langle \omega_\alpha \Omega_N^{-1} \omega_\beta \rangle
		\end{align*}
		
	\noindent
	and 
	
		\begin{align*}
			J_N(d) &= \frac{N^{2d}}{d!}  
				 \sum_{\substack{\lambda \vdash d \\ \ell(\lambda) \leq N}} s_\lambda(c_1,\dots,c_N) \frac{(\dim \mathsf{V}^\lambda)^2}{\dim \mathsf{W}^\lambda} \\
				&=  N^d\P(\LIS_d \leq N)\sum_{\beta \vdash d} p_\beta(c_1,\dots,c_N) \langle \Omega_N^{-1}\omega_\beta \rangle.
		\end{align*}

	\noindent
	The following bounds --- which say that $I_N(d)$ and $J_N(d)$ have maximum modulus in the case of trivial external fields --- will be needed in Section \ref{sec:Functional}.
		
	\begin{prop}
	\label{prop:BasicInequalities}
		For any $d,N \in \mathbb{N}$ and any $a_1,\dots,a_N,b_1,\dots,b_N,c_1,\dots,c_N \in \C$ of modulus
		at most one, we have
		
			\begin{equation*}
				|I_N(d)| \leq N^{2d} \quad\text{ and }\quad |J_N(d)| \leq \P(\LIS_d \leq N) N^{2d}.
			\end{equation*}
	\end{prop}
	
	\begin{proof}
		This follows from the fact that the Schur polynomials are monomial positive.
	\end{proof}

\section{Stable Asymptotics}
\label{sec:Stable}

In this section, we analyze the $N \to \infty$ asymptotics of the string coefficients of $I_N$ and $J_N.$ 
We obtain a convergent $N \to \infty$ asymptotic expansion for each fixed string coefficient, the coefficients of which 
count monotone walks on the symmetric groups with prescribed length and boundary conditions. 
These expansions are grouped together to form the stable HCIZ and BGW integrals $I$ and $J$, which are formal power series. 
The stable integrals $I$ and $J$ satisfy a formal power series version of Conjecture \ref{conj:Main},
the form of which points the way to an analytic solution.

	\subsection{String coefficients}
	Our present goal is to determine the $N \to \infty$ asymptotics of the string
	coefficients of $I_N$ and $J_N,$
	
		\begin{equation*} 
			\P(\LIS_d \leq N)\langle \omega_\alpha \Omega_N^{-1} \omega_\beta \rangle,
		\end{equation*}
		
	\noindent
	in the regime where $\alpha,\beta \vdash d$ are fixed and $N \to \infty.$ In this regime we may assume $N \geq d,$ so that
	the string coefficients are pure Plancherel expectations:

	\begin{equation*}
		\langle \omega_\alpha \Omega_N^{-1} \omega_\beta \rangle = \sum_{\lambda \vdash d} 
		\omega_\alpha(\lambda) \Omega_N^{-1}(\lambda) \omega_\beta(\lambda) \frac{(\dim \mathsf{V}^\lambda)^2}{d!}.
	\end{equation*}

	\noindent	
	Since 
	
		\begin{equation*}
			\lim_{N \to \infty} \Omega_N^{-1}(\lambda) = 1
		\end{equation*}
		
	\noindent
	for any fixed $\lambda,$ these Plancherel expectations are deformations of the usual inner product on the 
	center of the group algebra $\C\group{S}(d),$ with respect to which the functions $\omega_\alpha$ form an orthogonal basis:
	
		\begin{equation*}
			\langle \omega_\alpha \omega_\beta \rangle = \sum_{\lambda \vdash d} \omega_\alpha(\lambda)\omega_\beta(\lambda)
			\frac{(\dim \mathsf{V}^\lambda)^2}{d!} = \delta_{\alpha\beta} |C_\alpha|.
		\end{equation*}
		
	\noindent
	We thus have

		\begin{equation*}
			\langle \omega_\alpha \Omega_N^{-1} \omega_\beta \rangle = \delta_{\alpha\beta} |C_\alpha| + o(1)
		\end{equation*}
	
	\noindent
	as $N \to \infty,$ simply because of the orthogonality of irreducible characters. We will now quantify the error term
	in this approximation.
		
	Let $\hbar$ be a complex parameter, and consider the function on Young diagrams $\lambda$ defined by 
	
		\begin{equation*}
			\Psi_\hbar(\lambda) = \prod_{\Box \in \lambda} (1-\hbar c(\Box)).
		\end{equation*}
		
	\noindent
	This is a polynomial function of $\hbar$ whose roots are the reciprocals of the contents of the off-diagonal 
	cells of $\lambda.$ Explicitly, this polynomial is given by
		
		\begin{equation*}
			\Psi_\hbar(\lambda) = \sum_{r=0}^d (-\hbar)^r e_r(\lambda),
		\end{equation*}
		
	\noindent
	where $e_r(\lambda)$ is the degree $r$ elementary symmetric polynomial in $d$ variables,
	
		\begin{equation*}
			e_r = \sum_{\substack{ i \colon [r] \to [d] \\ i \text{ strictly increasing}}} x_{i(1)} \dots x_{i(r)},
		\end{equation*}
		
	\noindent
	evaluated on the contents of the diagram $\lambda$. Note that $e_r(\lambda)$ is a shifted symmetric
	function of $\lambda;$ see \cite{OP} for a discussion of shifted symmetric functions.	
	
	The function 
	
		\begin{equation*}
			\Psi_\hbar^{-1}(\lambda) = \frac{1}{\Psi_\hbar(\lambda)}
		\end{equation*}
		
	\noindent
	is a nonvanishing rational function of $\hbar$ whose poles are the roots of 
	$\Psi_\hbar(\lambda).$ In particular, for any diagram $\lambda,$ the function 
	$\Psi_\hbar^{-1}(\lambda)$ is analytic on the disc
	
		\begin{equation*}
			|\hbar| < \frac{1}{\max(\lambda_1-1,\ell(\lambda)-1)}
		\end{equation*}
		
	\noindent
	with Maclaurin series
	
		\begin{equation*}
			\Psi_\hbar^{-1}(\lambda) = \sum_{r=0}^\infty \hbar^r f_r(\lambda),
		\end{equation*}
		
	\noindent
	where $f_r(\lambda)$ is the degree $r$ complete symmetric polynomial in $d$ variables,
	
		\begin{equation*}
			f_r = \sum_{\substack{ i \colon [r] \to [d] \\ i \text{ weakly increasing}}} x_{i(1)} \dots x_{i(r)},
		\end{equation*}
		
	\noindent
	evaluated on the contents of $\lambda.$ 
	
	The functions $\Omega_N^{\pm 1}$ are recovered from 
	the functions $\Psi_\hbar^{\pm 1}$ by setting
	
		\begin{equation*}
			 \hbar = - \frac{1}{N}.
		\end{equation*}
		
	\noindent
	In particular, the $N \to \infty$ asymptotics of $\langle \omega_\alpha \Omega_N^{-1} \omega_\beta \rangle$
	may be obtained from the $\hbar \to 0$ asymptotics of $\langle \omega_\alpha \Psi_\hbar^{-1} \omega_\beta \rangle$,
	i.e. from its Taylor expansion around $\hbar = 0$ as derived above. We will now give a diagrammatic interpretation of this
	Maclaurin series.
		
	For a given pair of Young diagrams $\alpha,\beta \vdash d,$ we have the Taylor series
	
		\begin{equation*}
			\langle \omega_\alpha \Psi_\hbar^{-1} \omega_\beta \rangle = 
			 \sum_{r=0}^\infty \hbar^r \langle \omega_\alpha 
			f_r \omega_\beta \rangle,
		\end{equation*}
		
	\noindent
	which is absolutely convergent for $|\hbar| < \frac{1}{d-1}.$
	For any $\lambda \vdash d$, the observable $\omega_\alpha(\lambda)
	f_r(\lambda)\omega_\beta(\lambda)$ is the eigenvalue of the central element 
	$C_\alpha f_r(X_1,\dots,X_d) C_\beta$ acting in the irreducible representation $\mathsf{V}^\lambda$
	of the group algebra $\C\group{S}(d)$ corresponding to $\lambda$, where 
	$X_1,\dots,X_d \in \C\group{S}(d)$ are the Jucys-Murphy elements \cite{DG,OV}:
	
		\begin{equation*}
			X_j = \sum_{i < j} (i\ j), \quad 1 \leq j \leq d.
		\end{equation*}
		
	\noindent
	Thus, by the Fourier isomorphism,
	
		\begin{equation*}
			\C\group{S}(d) \simeq \bigoplus_{\lambda \vdash d} \mathrm{End} \mathsf{V}^\lambda,
		\end{equation*}
		
	\noindent
	the Plancherel expectation $\langle \omega_\alpha f_r \omega_\beta \rangle$ is the 
	normalized character of the central element $C_\alpha f_r(X_1,\dots,X_d) C_\beta$ in the 
	regular representation of $\C\group{S}(d),$
	i.e. the coefficient of the identity permutation in the sum
	
		\begin{equation*}
			\sum_{\rho \in C_\alpha, \sigma \in C_\beta}
			\sum_{\substack{ i,j \colon [r] \to [d] \\ i<j \text{ pointwise} \\ j \text{ weakly increasing}} } 
			\rho \left( i(1)\ j(1) \right) \dots \left( i(r)\ j(r) \right) \sigma.
		\end{equation*}
		
	\noindent
	Adopting the convention that permutations are multiplied from left to right, this number may be 
	visualized as follows. 
	
\begin{figure}
			\includegraphics{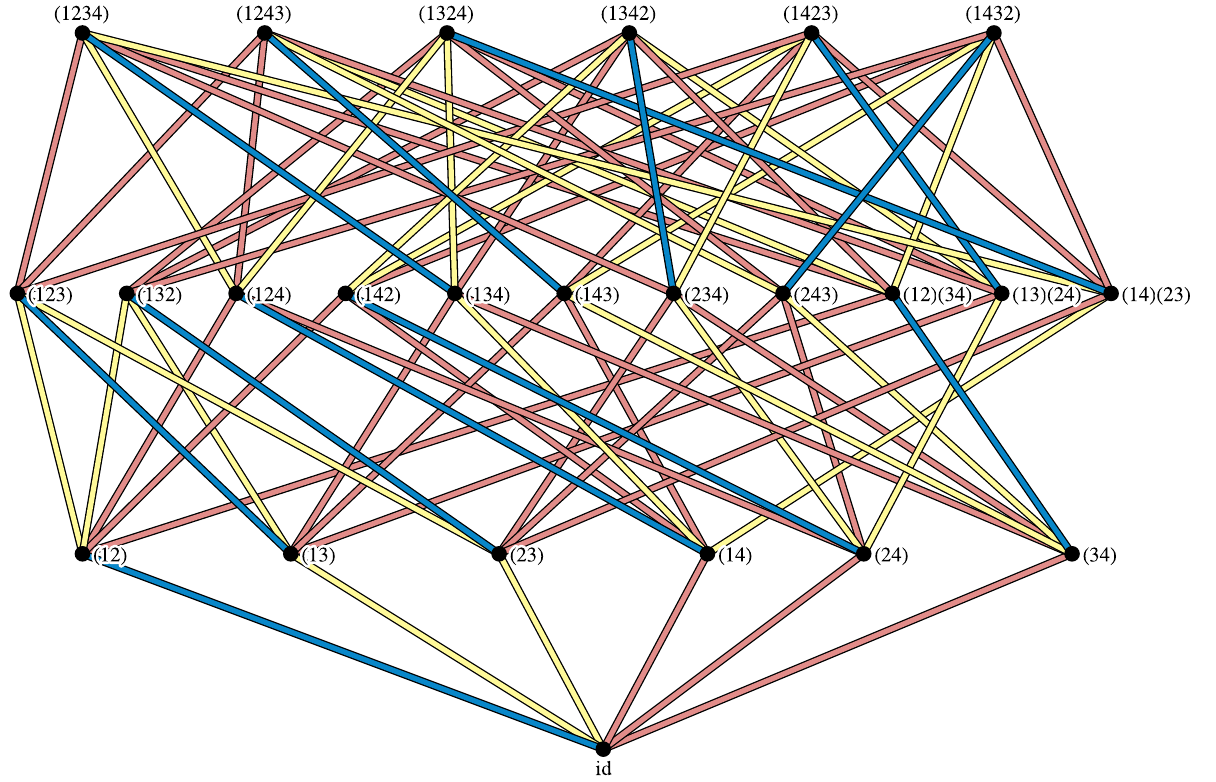}
			\caption{\label{fig:Cayley} Edge labeled Cayley graph of $\group{S}(4).$ Figure by M. LaCroix.}
		\end{figure}

	Identify the symmetric group $\group{S}(d)$ with its right Cayley graph, 
	as generated by the conjugacy class of transpositions. Introduce an edge labeling on this graph by marking each edge corresponding 
	to the transposition $(i\ j)$ with $j$, the larger of the two elements interchanged. 
	Thus, emanating from each vertex of $\group{S}(d)$, one sees a single 
	$2$-edge, two $3$-edges, three $4$-edges, etc. Figure \ref{fig:Cayley} shows $\group{S}(4)$ equipped 
	with this edge labeling. A walk on $\group{S}(d)$ is said to be \emph{monotone} if the labels of the edges it traverses form a
	weakly increasing sequence. Given Young diagrams $\alpha,\beta \vdash d$ and a nonnegative integer
	$r$, let $\vec{W}^r(\alpha,\beta)$ denote the number of $r$-step monotone walks on $\group{S}(d)$ which
	begin at a point of $C_\alpha$ and end at a point of $C_\beta$. Then from the calculation above
	we have the identity
	
		\begin{equation*}
			\langle \omega_\alpha f_r \omega_\beta \rangle = \vec{W}^r(\alpha,\beta).
		\end{equation*}
		
	\noindent
	Equivalently,
	
		\begin{equation*}
			\langle \omega_\alpha \Psi_\hbar^{-1} \omega_\beta \rangle =
			\sum_{r=0}^\infty \hbar^r \vec{W}^r(\alpha,\beta),
		\end{equation*}	
		
	\noindent
	the generating function for monotone walks on $\group{S}(d)$ with boundary conditions $\alpha,\beta,$
	the sum being absolutely convergent for $|\hbar| < \frac{1}{d-1}.$ In the special case $\alpha=(1^d),$ we
	have 
	
		\begin{equation*}
			\langle \Psi_\hbar^{-1} \omega_\beta \rangle =
			\sum_{r=0}^\infty \hbar^r \vec{W}^r(\beta),
		\end{equation*}
		
	\noindent
	where $\vec{W}^r(\beta) = \vec{W}^r(1^d,\beta)$ is the number of $r$-step monotone walks on $\group{S}(d)$
	which begin at the identity permutation and end at a permutation of cycle type $\beta.$ We may thus conclude
	the following.
	
		\begin{thm}
		\label{thm:FeynmanExpansions}
		For any positive integers $1 \leq d \leq N$ and any Young diagrams $\alpha,\beta \vdash d$,
		we have
	
			\begin{equation*}
				\langle \omega_\alpha \Omega_N^{-1} \omega_\beta \rangle 
				= \sum_{r=0}^\infty (-1)^r \frac{\vec{W}^r(\alpha,\beta)}{N^r} \quad\text{ and }\quad
				\langle \Omega_N^{-1} \omega_\beta \rangle 
				= \sum_{r=0}^\infty (-1)^r \frac{\vec{W}^r(\beta)}{N^r} 
			\end{equation*}
		
		\noindent
		and the series are absolutely convergent. 	
		\end{thm}
		
	Note that the $1/N$ expansions in Theorem \ref{thm:FeynmanExpansions} are not actually alternating series:
	their nonzero terms are either all negative or all positive.
		
	Theorem \ref{thm:FeynmanExpansions} gives a convergent $N \to \infty$ asymptotic expansion of 
	the Plancherel expectation $\langle \omega_\alpha \Omega_N^{-1} \omega_\beta \rangle$ wherein monotone
	walks play the role of Feynman diagrams. As a consistency check, observe that 
	
		\begin{equation*}
			\vec{W}^0(\alpha,\beta) = \delta_{\alpha\beta} |C_\alpha|,
		\end{equation*}
		
	\noindent
	corresponding to the fact that there exists a zero-step walk from $C_\alpha$ to $C_\beta$ if and only if these otherwise
	disjoint sets are equal, and in this case the number of such walks is just the cardinality of $C_\alpha.$ 
	
	\subsection{Stable integrals}
	For any fixed $N \in \mathbb{N},$ Theorem \ref{thm:FeynmanExpansions} describes the first $N$ nonconstant terms
	in the string expansions of $I_N$ and $J_N$: we have
		
		\begin{equation*}
			I_N = 1 + \sum_{d=1}^N \frac{z^d}{d!} \sum_{\alpha,\beta \vdash d} p_\alpha(a_1,\dots,a_N)
			p_\beta(b_1,\dots,b_N) \sum_{r=0}^\infty (-1)^r \frac{\vec{W}^r(\alpha,\beta)}{N^r} + 
			\text{ higher terms},
		\end{equation*}
		
	\noindent
	and 
	
		\begin{equation*}
			J_N = 1 + \sum_{d=1}^N \frac{z^{2d}}{d!} N^d\sum_{\beta \vdash d}
			p_\beta(b_1,\dots,b_N) \sum_{r=0}^\infty (-1)^r \frac{\vec{W}^r(\beta)}{N^r} + 
			\text{ higher terms}.
		\end{equation*}
		
	\noindent
	This description suggests that, as $N \to \infty,$ the integrals $I_N$ and $J_N$ approximate 
	generating functions for monotone walks on all of the symmetric groups, of all possible lengths and 
	boundary conditions. Unfortunately, for any finite $N,$ almost all terms of the string expansion are ``higher terms.'' 
	
	To see past this analytic limitation, let us view $z$ as a formal variable, and replace the number 
	$-\frac{1}{N}$ with a formal semiclassical parameter $\hbar.$ Furthermore, let us replace the eigenvalues of the matrices
	$A,B,$ and $C=AB$ with countably infinite alphabets of commuting indeterminates, these being formal stand-ins for the eigenvalues 
	infinite-dimensional matrices. Let $\Lambda_A,\Lambda_B,\Lambda_C$ be the affiliated algebras of symmetric functions,
	i.e. the polynomial algebras
	
		\begin{equation*}
			\Lambda_A = \C[p_1(A),p_2(A),\dots],\quad \Lambda_B = \C[p_1(B),p_2(B),\dots],\quad \Lambda_C = \C[p_1(C),p_2(C),\dots],
		\end{equation*}
		
	\noindent
	where
	
		\begin{equation*}
			p_k(A) = \sum_{a \in A} a^k,\quad p_k(B) = \sum_{b \in B} b^k,\quad p_k(C) = \sum_{c \in C} c^k, \qquad k \in \mathbb{N},
		\end{equation*}
		
	\noindent
	are the pure power sums over these alphabets. Set $\Lambda_{A,B} = \Lambda_A \otimes \Lambda_B.$ 
	
	We define the stable HCIZ integral to be the formal power series  
		
		\begin{equation*}
			I = 1 + \sum_{d=1}^\infty \frac{z^d}{d!} \sum_{\alpha,\beta \vdash d} p_\alpha(A) p_\beta(B)
			\sum_{r=0}^\infty \hbar^r \vec{W}^r(\alpha,\beta),
		\end{equation*}
		
	\noindent
	which is an element of the ring $\Lambda_{A,B}[[z,\hbar]].$
	Similarly, we define the stable BGW integral to be the formal power series
	
		\begin{equation*}
			J  = 1 + \sum_{d=1}^\infty \frac{z^{2d}}{d!} (-1)^d\hbar^{-d} \sum_{\beta \vdash d}p_\beta(B)
			\sum_{r=0}^\infty \hbar^r \vec{W}^r(\beta),
		\end{equation*}	
		
	\noindent
	which is an element of $\Lambda_C[[z,\hbar^{\pm 1}]].$
	Thus $I$ and $J$ are ``grand canonical'' partition functions enumerating monotone walks of all possible lengths and boundary conditions, over all
	symmetric groups.
	
		\begin{thm}
		\label{thm:StableExponentials}
			We have 
			
			\begin{equation*}
				I = e^F \quad\text{ and }\quad J=e^G,
			\end{equation*}
			
		\noindent
		where
		
			\begin{equation*}
				F = \sum_{d=1}^\infty \frac{z^d}{d!} \sum_{\alpha,\beta \vdash d} 
				p_\alpha(A) p_\beta(B) \sum_{r=0}^\infty \hbar^r \mon^r(\alpha,\beta)
			\end{equation*}
		
		\noindent	
		and $\mon^r(\alpha,\beta)$ is the number of monotone $r$-step walks on $\group{S}(d)$ which begin
		at a permutation of cycle type $\alpha,$ end at a permutation of cycle type $\beta,$ and have the property
		that their steps and endpoints together generate a transitive subgroup of $\group{S}(d),$ and
		
			\begin{equation*}
				G = \sum_{d=1}^\infty \frac{z^d}{d!} (-1)^d \hbar^{-d}\sum_{\beta \vdash d} 
				p_\beta(C) \sum_{r=0}^\infty \hbar^r \mon^r(\beta)
			\end{equation*}
			
		\noindent
		with $\mon^r(\beta) = \mon^r(1^d,\beta).$
		\end{thm}
		
	Theorem \ref{thm:StableExponentials} follows from a fundamental result in enumerative combinatorics,
	the Exponential Formula \cite[Chapter 5]{Stanley:EC2}, according to which the exponential of a 
	generating function for a class of ``connected'' combinatorial structures is a generating function
	for possibly ``disconnected'' structures of the same type. For walks on groups, the role of connectedness is played by transitivity.
	For a careful justification of the use of the Exponential Formula in the context of monotone walks on symmetric
	groups, see \cite{GGN1,GGN2,GGN4}.
	
	\subsection{Topological expansion}
	The numbers $\mon^r(\alpha,\beta)$ and $\mon^r(\beta)$ appearing in Theorem \ref{thm:StableExponentials}
	are known as the monotone double and single Hurwitz numbers, respectively. These enumerative quantities,
	introduced in \cite{GGN1,GGN2,GGN3} and studied in numerous articles since, are a combinatorial variant
	of the classical double and single Hurwitz numbers $H^r(\alpha,\beta)$ and $H^r(\beta)=H^r(1^d,\beta),$ which count transitive
	$r$-step walks $C_\alpha \to C_\beta$ without the monotonicity constraint. Clearly, $\mon^r(\alpha,\beta) \leq H^r(\alpha,\beta),$
	and in a sense monotone Hurwitz numbers are a ``desymmetrized'' version of classical Hurwitz numbers;
	see \cite{GGN1,GGN2}.
	
	Reversing a classical construction due to Hurwitz \cite{Hurwitz} (see \cite{EEHS} for a modern treatment), 
	the number $H^r(\alpha,\beta)$ may alternatively be interpreted as the number of isomorphism
	classes of degree $d$ branched covers of the Riemann sphere $\mathbf{P}^1(\C)$ by a compact, connected
	Riemann surface $\mathbf{S}$ which have profiles $\alpha,\beta \vdash d$
	over $0,\infty \in \mathbf{P}^1(\C)$ and the simplest nontrivial branching over the $r$th roots of unity on the sphere.
	The monotone double Hurwitz number $\mon^r(\alpha,\beta)$ is a signed enumeration of the same class of covers, see
	\cite{ACEH,MN}. The genus $g$ of $\mathbf{S}$ is determined by the data $d,r,\alpha,\beta$ according to the Riemann-Hurwitz formula,
	
		\begin{equation*}
			g = \frac{r+2-\ell(\alpha)-\ell(\beta)}{2},
		\end{equation*}
		
	\noindent
	with the understanding that $H^r(\alpha,\beta)=0$ unless this formula returns a nonnegative integer. In particular, one may 
	parameterize nonzero (classical and monotone) Hurwitz numbers by genus, setting $H_g(\alpha,\beta) := H^{2g-2+\ell(\alpha)+\ell(\beta)}(\alpha,\beta)$
	and $\mon_g(\alpha,\beta): = \mon^{2g-2+\ell(\alpha)+\ell(\beta)}(\alpha,\beta).$ In the genus parameterization,
	Theorem \ref{thm:StableExponentials} becomes the following topological expansion of the stable 
	HCIZ and BGW integrals.

	\begin{thm}
		\label{thm:StableTopExpansions}
		We have 
			
			\begin{equation*}
				I = e^{\sum_{g=0}^\infty \hbar^{2g-2} F^{(g)}} \quad\text{ and }\quad
				J = e^{\sum_{g=0}^\infty \hbar^{2g-2} G^{(g)}},
			\end{equation*}
			
		\noindent
		where
		
			\begin{equation*}
				F^{(g)} = \sum_{d=1}^\infty \frac{z^d}{d!} \sum_{\alpha,\beta \vdash d} 
				\hbar^{\ell(\alpha)+\ell(\beta)} p_\alpha(A) p_\beta(B) \mon_g(\alpha,\beta).
			\end{equation*}
		
		\noindent
		and
		
			\begin{equation*}
				G^{(g)} = \sum_{d=1}^\infty \frac{z^{2d}}{d!} (-1)^d \sum_{\beta \vdash d} 
				\hbar^{\ell(\beta)} p_\beta(C) \mon_g(\beta).
			\end{equation*}
	\end{thm}
	
	\begin{proof}
		Applying the Riemann-Hurwitz formula, we have
		
			\begin{align*}
				F &= \sum_{d=1}^\infty \frac{z^d}{d!} \sum_{\alpha,\beta \vdash d} p_\alpha(A) p_\beta(B) \sum_{g=0}^\infty \hbar^{2g-2+\ell(\alpha)+\ell(\beta)}\mon_g(\alpha,\beta) \\
				&= \sum_{g=0}^\infty \hbar^{2g-2} \sum_{d=1}^\infty \frac{z^d}{d!} \sum_{\alpha,\beta \vdash d} \hbar^{\ell(\alpha)+\ell(\beta)}p_\alpha(A) p_\beta(B) \mon_g(\alpha,\beta)
			\end{align*}	
			
		\noindent
		and 
		
			\begin{align*}
				G &= \sum_{d=1}^\infty \frac{z^d}{d!} (-1)^d \hbar^{-d}\sum_{\beta \vdash d} 
				p_\beta(B) \sum_{g=0}^\infty \hbar^{2g-2+d+\ell(\beta)} \mon^r(\beta) \\
				&= \sum_{g=0}^\infty \hbar^{2g-2} \sum_{d=1}^\infty \frac{z^d}{d!} (-1)^d 
				\sum_{\beta \vdash d} \hbar^{\ell(\beta)} p_\beta(C) \mon_g(\beta).
			\end{align*}
	\end{proof}	
	
	\subsection{Topological factorization}
	For any nonnegative integer $k,$ the topological expansions of $I$ and $J$ given by 
	Theorem \ref{thm:StableTopExpansions} can be split into two corresponding factors,
	
		\begin{equation*}
			I = e^{\sum_{g=0}^k \hbar^{2g-2} F^{(g)}} e^{\sum_{g=k+1}^\infty \hbar^{2g-2} F^{(g)}} 
			\quad\text{ and }\quad
			J = e^{\sum_{g=0}^k \hbar^{2g-2} G^{(g)}} e^{\sum_{g=k+1}^\infty \hbar^{2g-2} G^{(g)}}.
		\end{equation*}
		
	\noindent
	These factorizations have a clear enumerative meaning: the first factor is a generating fucntion enumerating 
	possibly disconnected covers/walks in which each connected component has genus at most $k$, while the second factor
	is a generating function enumerating possibly disconnected covers/walks in which each connected component has genus at least $k+1.$ 
	This may be equivalently stated as follows. Define the disconnected monotone double and single
	Hurwitz numbers by 
	
		\begin{equation*}
			\mon_g^\bullet(\alpha,\beta) = \vec{W}^{r_g(\alpha,\beta)}(\alpha,\beta) 
			\quad\text{ and }\quad
			\mon_g^\bullet(\beta) = \mon_g^\bullet(1^d,\beta),
		\end{equation*}
		
	\noindent
	where $g \in \mathbb{Z}$ and $r_g(\alpha,\beta) =2g-2+\ell(\alpha)+\ell(\beta).$ In particular,
	for disconnected Hurwitz numbers the genus $g$ may be negative (corresponding to the fact
	that the Euler characteristic is additive), but $\mon_g^\bullet(\alpha,\beta)$ vanishes unless
	$r_g(\alpha,\beta) \geq 0.$ In terms of disconnected monotone Hurwitz numbers, the above
	factorization identities may be stated as follows.

	\begin{thm}
	\label{thm:StableTopFactorization}
	For any $k \in \mathbb{N} \cup \{0\}$	
		\begin{equation*}
			\frac{I}{e^{\sum_{g=0}^k \hbar^{2g-2} F^{(g)}}} = 1 + \sum_{d=1}^\infty \frac{z^d}{d!} \sum_{\alpha,\beta \vdash d} 
			\hbar^{\ell(\alpha)+\ell(\beta)} p_\alpha(A)p_\beta(B) \sum_{g=k+1}^\infty
			\hbar^{2g-2} \mon_g^\bullet(\alpha,\beta)
		\end{equation*}
		
	\noindent
	and
	
		\begin{equation*}
			\frac{J}{e^{\sum_{g=0}^k \hbar^{2g-2} G^{(g)}}} = 1 + \sum_{d=1}^\infty \frac{z^{2d}}{d!} \sum_{\alpha,\beta \vdash d} (-1)^d
			\hbar^{\ell(\beta)}p_\beta(B) \sum_{g=k+1}^\infty
			\hbar^{2g-2} \mon_g^\bullet(\beta).
		\end{equation*}
		
	\end{thm}
	
	As a corollary of Theorem \ref{thm:StableTopFactorization}, we obtain the following pair 
	of ``topological bounds,'' which are significant since they indicate what sorts of bounds we may expect to be valid for
	the entire functions $I_N$ and $J_N$ at finite $N.$ Let us introduce the following formal order notation: given a formal power
	series $Z \in \Lambda_{A,B}[[z,\hbar]]$ and a nonnegative integer $n,$ we write 
	
		\begin{equation*}
			Z = O(\hbar^n)
		\end{equation*}
		
	\noindent
	if $Z$ belongs to the principal ideal generated by $\hbar^n.$ We use the analogous 
	order notation in $\Lambda_C[[z,\hbar^{\pm 1}]].$
	
		\begin{cor}
		\label{cor:StableTopBound}
		For each $k \in \mathbb{N} \cup \{0\},$ 
		
			\begin{equation*}
				1 - \frac{I}{e^{\sum_{g=0}^k \hbar^{2g-2} F^{(g)}}} = O(\hbar^{2k}) \quad\text{ and }\quad
				1 - \frac{J}{e^{\sum_{g=0}^k \hbar^{2g-2} G^{(g)}}} = O(\hbar^{2k}).
			\end{equation*}
		\end{cor}	
	
\section{Functional Asymptotics}
\label{sec:Functional}
In this Section, we prove our main result, Theorem \ref{thm:Main}. To achieve this,
we must bridge the gap between the $N<\infty$ string expansions of the HCIZ and BGW
integrals and their $N=\infty$ stable topological expansions. It is here that the mollifying effect
of the $\LIS$ distribution plays a critical role: it controls the tail of the finite $N$ 
string expansions of $I_N$ and $J_N,$ effectively truncating them to polynomials 
of degree $O(N^2)$ for small $z$. The existence of this quadratic cutoff is a key
feature of $I_N$ and $J_N$ that has not previously been recognized.

\subsection{Analytic candidates}
Throughout this section, we will use the following notation for complex polydiscs. 
Given a real number $\rho$ and a positive integer $N$, we will ambiguously write
$\overline{\D}_\rho^N$ to mean either of the closed polydiscs

	\begin{equation*}
		\overline{\D}_\rho \times \overline{\D}_1^N \times \overline{\D}_1^N
		\quad\text{ or }\quad
		\overline{\D}_\rho \times \overline{\D}_1^N,
	\end{equation*}
	
\noindent
where $\overline{\D}_\rho$ is the closed origin-centred disc of radius $\rho$ in the complex plane. Although 
the first of these domains lives in $\mathbb{C}^{2N+1}$ and 
the second lives in $\mathbb{C}^{N+1},$
which of the two domains $\overline{\D}_\rho^N$ is intended to represent will be clear from context. 
Let $\|\cdot\|_\rho$ denote the sup norm on $\overline{\D}_\rho^N$. Note that this is really a sequence of 
norms defined on a sequence of domains of growing dimension. 

Let $N \in \mathbb{N}$ be a positive integer, and let $a_1,\dots,a_N,b_1,\dots,b_N,c_1,\dots,c_N$ be any points sampled
from the closed unit disc in $\C.$ Consider the corresponding specializations
	
		\begin{equation*}
			\Lambda_{A,B}[[z,\hbar]] \to \C[[z]] \quad\text{ and }\quad \Lambda_C[[z,\hbar]] \to \C[[z]]
		\end{equation*}
		
	\noindent
	defined by setting $\hbar=-1/N$ and
	
		\begin{equation*}
			A=\{a_1,\dots,a_N\}, B=\{b_1,\dots,b_N\}, C=\{c_1,\dots,c_N\},
		\end{equation*}
		
	\noindent
	and let

	\begin{align*}
		F_N^{(g)} &= \sum_{d=1}^\infty \frac{z^d}{d!}
				  \sum_{\alpha,\beta \vdash d} 
				\frac{p_\alpha(a_1,\dots,a_N)}{N^{\ell(\alpha)}}
				\frac{p_\beta(b_1,\dots,b_N)}{N^{\ell(\beta)}} (-1)^{\ell(\alpha)+\ell(\beta)} \mon_g(\alpha,\beta),  \\
		G_N^{(g)} &= \sum_{d=1}^\infty \frac{z^{2d}}{d!}
				\sum_{\beta \vdash d} 
				 \frac{p_\beta(c_1,\dots,c_N)}{N^{\ell(\beta)}} (-1)^{d+\ell(\beta)} \mon_g(\beta).
	\end{align*}
	
\noindent
be the images of $F^{(g)}$ and $G^{(g)}$ under these specializations.
A priori, $F_N^{(g)}$ and $G_N^{(g)}$ are only formal power series. However, they are in fact absolutely summable,
and hence define analytic functions. This follows from an established result on the convergence of generating functions for monotone
Hurwitz numbers \cite{GGN5}.

	\begin{thm}
	\label{thm:GenusSpecific}
	For each $g \in \mathbb{N} \cup \{0\},$ the power series 
	
		\begin{align*}
			\mon_g^{\text{simple}} &= \sum_{d=1}^\infty \frac{z^d}{d!} \mon_g(1^d,1^d), \\
			\mon_g^{\text{single}} &= \sum_{d=1}^\infty \frac{z^d}{d!} \sum_{\beta \vdash d} \mon_g(1^d,\beta), \\
			\mon_g^{\text{double}} &= \sum_{d=1}^\infty \frac{z^d}{d!} \sum_{\alpha,\beta \vdash d} \mon_g(\alpha,\beta) 
		\end{align*}
		
	\noindent
	have radii of convergence exactly $2/27,$ at least $1/27,$ and at least $1/54,$ respectively.
	\end{thm}
	
\noindent
The exact computation of the radius of convergence of the generating function for monotone simple Hurwitz 
numbers follows from a rational parameterization of this series in terms of the Gauss hypergeometric function, 
see \cite{GGN2,GGN5}. A combinatorial argument based on sorting transpositions (a variant of the Hurwitz
braid action) then shows that the radius
of convergence drops by at most a factor of two for each new branch point added, see \cite{GGN5} for details.
The author believes that the radius of convergence is in fact exactly $2/27$ in all three cases, but this has not
been proved.

It was pointed out to the author by Philippe Di Francesco that the number of isomorphism classes of finite groups 
of order $p^N,$ with $p$ prime, is known \cite{Pyber}  to be asymptotically $p^{\frac{2}{27}N^3}$ as $N \to \infty.$
The author has no explanation for this numerical coincidence. For another interesting appearance of the number
$2/27,$ see \cite{KT}.

	\begin{thm}
	\label{thm:AnalyticCandidates}
		There exists $\delta > 0$ such that the series 
		$F_N^{(g)}$ and $G_N^{(g)}$ converge absolutely and uniformly on $\overline{\D}_\delta^N$,
		for all $g \geq 0$ and $N \geq 1.$
	\end{thm}
	
	\begin{proof}
		For any Young diagrams $\alpha,\beta$, we have
		
			\begin{equation*}
				\left| \frac{p_\alpha(a_1,\dots,a_N)}{N^{\ell(\alpha)}} \right|, 
				\left| \frac{p_\beta(b_1,\dots,b_N)}{N^{\ell(\beta)}} \right|,
				\left| \frac{p_\beta(c_1,\dots,c_N)}{N^{\ell(\beta)}} \right| \leq 1,
			\end{equation*} 
			
		\noindent
		on $\overline{\D}_\delta^N$. We thus have
		
			\begin{align*}
				|F_N^{(g)}| &\leq \sum_{d=1}^n \frac{\delta^d}{d!} \sum_{\alpha,\beta \vdash d} \mon_g(\alpha,\beta) \\
				|G_N^{(g)}| &\leq \sum_{d=1}^n \frac{\delta^{2d}}{d!} \sum_{\beta \vdash d} \mon_g(\beta) 
			\end{align*}
		
		\noindent
		uniformly on $\overline{\D}_\delta^N$
		for any $n \in \mathbb{N},$ and the claim thus follows from Theorem \ref{thm:GenusSpecific}.
	\end{proof}

Let us fix $\delta >0$ so that $F_N^{(g)}$ and $G_N^{(g)}$ converge to define analytic functions on
on $\overline{\D}_\delta^N,$ for all $N \in \mathbb{N}.$ Then, these functions are uniformly bounded in the following sense.

	\begin{cor}
	\label{cor:UniformFreeEnergies}
	We have 
		
		\begin{equation*}
			\sup_{N \in \mathbb{N}} \|F_N^{(g)}\|_\delta < \infty
			\quad\text{ and }\quad
			\sup_{N \in \mathbb{N}} \|G_N^{(g)}\|_\delta  < \infty.
		\end{equation*}
	\end{cor}
		
\subsection{Polynomial approximation}
We now consider the behavior of the full string expansions of $I_N$ and $J_N$ given by Theorem \ref{thm:StringExpansions}
with $N$ large but finite. In this regime, the factor $\P(\LIS_d \leq N)$ has a dramatic effect on the 
string expansions --- it effectively truncates them to polynomials of degree $O(N^2)$. The mechanism behind this cutoff
is the law of large numbers for longest increasing subsequences in random permutations, 
which is due to Vershik and Kerov \cite{Kerov}: we have

	\begin{equation*}
		\label{eqn:LISLLN}
		\lim_{d \to \infty} \frac{\LIS_d}{\sqrt{d}} = 2, 
	\end{equation*}
	
\noindent
where the convergence is in probability. A detailed exposition of this LLN is given in \cite{Romik}, 
which also presents the corresponding central limit theorem of Baik-Deift-Johansson \cite{BDJ}, which
asserts Tracy-Widom fluctuations of $\LIS_d$ around $2\sqrt{d}$ on the scale $d^{1/6}$. In particular, 
the distribution of the longest increasing subsequence in large uniformly random permutation is strongly
concentrated around its mean. This implies that, for large $N$ we have the approximate step function behavior

	\begin{equation*}
		\P(\LIS_d \leq N) \approx \begin{cases}
			1, \quad 1 \leq d \leq \frac{1}{4}N^2 \\
			0, \quad d > \frac{1}{4}N^2
		\end{cases}.
	\end{equation*}
	
\noindent
Consequently, for $|z|$ small and $N$ large, $I_N$ and $J_N$ are well-approximated by their
``string polynomials''

	\begin{equation*}
		\tilde{I}_N = 1 + \sum_{d=1}^{\lfloor \frac{1}{4}N^2 \rfloor} \frac{z^d}{d!} \sum_{\alpha,\beta \vdash d} 
		p_\alpha(a_1,\dots,a_N) p_\beta(b_1,\dots,b_N) \langle \omega_\alpha \Omega_N^{-1} \omega_\beta \rangle 
	\end{equation*}
	
\noindent
and 

	\begin{equation*}
		\tilde{J}_N = 1 + \sum_{d=1}^{\lfloor \frac{1}{4}N^2 \rfloor} \frac{z^{2d}}{d!} N^d\sum_{\beta \vdash d} 
		p_\beta(c_1,\dots,c_N) \langle \Omega_N^{-1} \omega_\beta \rangle,
	\end{equation*}
	
\noindent
which are obtained from the string expansions of $I_N$ and $J_N$ as given by Theorem \ref{thm:StringExpansions}
by replacing the factor $\P(\LIS_d \leq N)$ with the above step function.
		
\subsection{Feynman extension}
The polynomial approximations $\tilde{I}_N$ and $\tilde{J}_N$ of $I_N$ and $J_N$ are only useful insofar 
as we are able to understand the Plancherel expectations 

	\begin{equation*}
		\langle \omega_\alpha \Omega_N^{-1} \omega_\beta \rangle = \sum_{\substack{\lambda \vdash d \\ \ell(\lambda) \leq N}}
		\omega_\alpha(\lambda) \Omega_N^{-1}(\lambda) \omega_\beta(\lambda) \frac{(\dim \mathsf{V}^\lambda)^2}{|\group{S}_N(d)|}
	\end{equation*}
	
\noindent
in the range $1 \leq d \leq \frac{1}{4}N^2.$ This means that we must extend Theorem \ref{thm:FeynmanExpansions},
which gives the convergent $1/N$ expansion

	\begin{equation*}
		\langle \omega_\alpha \Omega_N^{-1} \omega_\beta \rangle = \sum_{r=0}^\infty \left( -\frac{1}{N} \right)^r \vec{W}^r(\alpha,\beta)
		= \frac{(-1)^{\ell(\alpha)+\ell(\beta)}}{N^{\ell(\alpha)+\ell(\beta)}} \sum_{\substack{g = -\infty \\ 2-2g \leq \ell(\alpha)+\ell(\beta)}}
		N^{2-2g} \mon_g^\bullet(\alpha,\beta)
	\end{equation*}
	
\noindent
in the linear range $1 \leq d \leq N,$ to the range where $d$ may be as large as $\frac{1}{4}N^2.$ This may be done as follows.

For any $d,N \in \mathbb{N}$ we may rewrite the expectation $\langle \omega_\alpha \Omega_N^{-1} \omega_\beta \rangle$
as a conditional expectation against the unrestricted Plancherel measure: we have

	\begin{equation*}
		\langle \omega_\alpha \Omega_N^{-1} \omega_\beta \rangle = \frac{1}{\P(\LIS_d \leq N)} \sum_{\substack{\lambda \vdash d \\ \ell(\lambda) \leq N}}
		\omega_\alpha(\lambda) \Omega_N^{-1}(\lambda) \omega_\beta(\lambda) \frac{(\dim \mathsf{V}^\lambda)^2}{d!}
	\end{equation*}
	
\noindent
Let us split this conditional expectation into two pieces: we write

	\begin{equation*}
		\langle \omega_\alpha \Omega_N^{-1} \omega_\beta \rangle = \langle \omega_\alpha \Omega_N^{-1} \omega_\beta \rangle_1 + 
		\langle \omega_\alpha \Omega_N^{-1} \omega_\beta \rangle_2,
	\end{equation*}
	
\noindent
where

	\begin{equation*}
		 \langle \omega_\alpha \Omega_N^{-1} \omega_\beta \rangle_1= \frac{1}{\P(\LIS_d \leq N)} \sum_{\substack{\lambda \vdash d \\ \ell(\lambda) \leq N \\ \lambda_1 \leq N}}
		\omega_\alpha(\lambda) \Omega_N^{-1}(\lambda) \omega_\beta(\lambda) \frac{(\dim \mathsf{V}^\lambda)^2}{d!}
	\end{equation*}
	
\noindent
and

	\begin{equation*}
		 \langle \omega_\alpha \Omega_N^{-1} \omega_\beta \rangle_2 =
		\frac{1}{\P(\LIS_d \leq N)} \sum_{\substack{\lambda \vdash d \\ \ell(\lambda) \leq N \\ \lambda_1 > N}}
		\omega_\alpha(\lambda) \Omega_N^{-1}(\lambda) \omega_\beta(\lambda) \frac{(\dim \mathsf{V}^\lambda)^2}{d!}.
	\end{equation*}
	
\noindent
In the range $1 \leq d \leq N,$ the second component of this decomposition vanishes. 
In the extended range $N < d \leq \frac{1}{4}N^2,$ when $N$ is large, 
the first component of this decomposition is virtually equal to $\langle \omega_\alpha \Omega_N^{-1} \omega_\beta \rangle,$ while
the second is negligible. Indeed, it follows from the Vershik-Kerov limit shape theorem \cite{Kerov,Romik} that for $N < d_N \leq
\frac{1}{4}N^2,$ a Plancherel-random Young diagram with $d_N$ cells is contained in the $N \times N$ rectangular diagram $R(N,N)$
with overwhelming probability.

Observe now that the massive component $\langle \omega_\alpha \Omega_N \omega_\beta \rangle_1$ of 
$\langle \omega_\alpha \Omega_N \omega_\beta \rangle$ admits an absolutely convergent $1/N$ expansion.
Indeed, for any $\lambda \subseteq R(N,N),$ we have the absolutely convergent expansion

	\begin{equation*}
		\Omega_N^{-1}(\lambda) = \prod_{\Box \in \lambda} \frac{1}{1 + \frac{c(\Box)}{N}} = \sum_{r=0}^\infty 
		\left( -\frac{1}{N} \right)^r f_r(\lambda),
	\end{equation*}
	
	\noindent
	so that 

		\begin{equation*}
			\langle \omega_\alpha \Omega_N^{-1} \omega_\beta \rangle_1
			= \sum_{r=0}^\infty \left( -\frac{1}{N} \right)^r \vec{W}^r_N(\alpha,\beta),
		\end{equation*}
	
	\noindent
	where 

		\begin{equation*}
			\vec{W}^r_N(\alpha,\beta) =
			\sum_{\substack{\lambda \vdash d \\ \lambda \subseteq R(N,N)}}
			\omega_\alpha(\lambda) f_r(\lambda) \omega_\beta(\lambda) \frac{(\dim \mathsf{V}^\lambda)^2}{d!}.
		\end{equation*}
	
	\noindent
	agrees with $\vec{W}^r(\alpha,\beta)$ up to an exponentially small error.
	Thus, for any fixed but arbitrary $s \in \mathbb{N} \cup \{0\},$ 
	we can replace the first $s$ coefficients of the massive component
	$\langle \omega_\alpha \Omega_N \omega_\beta \rangle_1$ with their 
	$N$-independent counterparts up to an exponentially small error. Ignoring
	the negligible component $\langle \omega_\alpha \Omega_N \omega_\beta \rangle_2,$
	this gives the $N \to \infty$ asymptotic approximation 
	
		\begin{equation*}
			\langle \omega_\alpha \Omega_N^{-1} \omega_\beta \rangle
			= \sum_{r=0}^s \left( -\frac{1}{N} \right)^r \vec{W}^r(\alpha,\beta) + O\left(\frac{1}{N^{s+1}}\right),
		\end{equation*}
		
	\noindent
	which extends Theorem \ref{thm:FeynmanExpansions} to the range $1 \leq d \leq \frac{1}{4}N^2.$
	Note that this expansion implies the sharper estimate
	
		\begin{equation*}
			\langle \omega_\alpha \Omega_N^{-1} \omega_\beta \rangle
			= \sum_{r=0}^s \left( -\frac{1}{N} \right)^r \vec{W}^r(\alpha,\beta) + O\left(\frac{1}{N^{s+2}}\right),
		\end{equation*}
		
	\noindent
	since the numbers $\vec{W}^r(\alpha,\beta)$ which are nonzero correspond to either $r$ even,
	or $r$ odd. In particular, for any $k \in \mathbb{N} \cup \{0\},$ we have that 
	
		\begin{equation*}
			\langle \omega_\alpha \Omega_N^{-1} \omega_\beta \rangle
			= \frac{(-1)^{\ell(\alpha)+\ell(\beta)}}{N^{\ell(\alpha)+\ell(\beta)}}
			\sum_{\substack{g = -\infty \\ 2-2g \leq \ell(\alpha) + \ell(\beta)}}^k 
			N^{2-2g} \mon_g^{\bullet}(\alpha,\beta) + O\left( N^{-2k}\right).
		\end{equation*}

\subsection{Topological bound}
The following topological bound bridges the gap between formal asymptotics and functional asymptotics.
Conceptually, this result is the unstable analytic shadow of the 
stable topological bounds appearing in Corollary \ref{cor:StableTopBound}. 

	\begin{thm}
	\label{thm:TopBound}
		There exists $\gamma > 0$ such that, for each fixed $k \in \mathbb{N} \cup \{0\},$ we have
						
				\begin{equation*}
					\bigg{\|} 1 - \frac{I_N}{e^{\sum_{g=0}^k N^{2-2g}F_N^{(g)}}} \bigg{\|}_\gamma = O(N^{2-2k})
					\quad\text{ and }\quad
					\bigg{\|} 1 - \frac{J_N}{e^{\sum_{g=0}^k N^{2-2g}G_N^{(g)}}} \bigg{\|} = O(N^{2-2k})
				\end{equation*}
				
		\noindent
		as $N \to \infty.$ 
	\end{thm}
	
	\begin{proof}	
		We give the proof for the HCIZ integral; the argument for the BGW integral is essentially the same.
		
		With $\delta$ as in Theorem \ref{thm:AnalyticCandidates}, take 
		$\gamma \leq \delta$ sufficiently small so that $I_N$ can be replaced with $\tilde{I}_N$ as 
		$N \to \infty.$ Replacing the coefficients of $\tilde{I}_N$ with their asymptotic expansions to 
		order $k+1$ and applying 	Theorem \ref{thm:StableTopFactorization}, we obtain
		
			\begin{align*}
				1-\frac{\tilde{I}_N}{e^{\sum_{g=0}^k N^{2-2g}F_N^{(g)}}} &= \sum_{d=1}^{\lfloor \frac{1}{4}N^2 \rfloor} 
				\frac{z^d}{d!} \sum_{\alpha,\beta \vdash d} \frac{p_\alpha(a_1,\dots,a_N)}{N^{\ell(\alpha)}}\frac{p_\beta(b_1,\dots,b_N)}{N^{\ell(\beta)}}
				\left(N^{-2k}\mon_{k+1}^\bullet(\alpha,\beta) + O(N^{-2k-2}) \right) \\
				&+ O(z^{\lfloor \frac{1}{4}N^2\rfloor + 1}).
			\end{align*}
			
		\noindent
		On $\overline{\D}_\gamma^N,$ we have the estimate 
		
			\begin{align*}
				&\bigg{|}\sum_{d=1}^{\lfloor \frac{1}{4}N^2 \rfloor} 
				\frac{z^d}{d!} \sum_{\alpha,\beta \vdash d} \frac{p_\alpha(a_1,\dots,a_N)}{N^{\ell(\alpha)}}\frac{p_\beta(b_1,\dots,b_N)}{N^{\ell(\beta)}}
				\left(N^{-2k}\mon_{k+1}^\bullet(\alpha,\beta) + O(N^{-2k-2}) \right) \bigg{|} \\
				&\leq N^{-2k} \sum_{d=1}^{\lfloor \frac{1}{4}N^2 \rfloor} 
				\frac{\gamma^d}{d!} \sum_{\alpha,\beta \vdash d} \left(\mon_{k+1}^\bullet(\alpha,\beta) + O(N^{-2})\right) \\
				&=O(N^{2-2k}).
			\end{align*}
			
		\end{proof}
		
	\begin{cor}
	\label{cor:NonVanishing}
		For $N$ sufficiently large, the integrals $l_N$ and $J_N$ are non-vanishing on $\overline{\D}_{\gamma}^N.$
		In particular, for $N$ sufficiently large, 
		$\log I_N$ and $\log J_N$ are defined and analytic on $\overline{\D}_{\gamma}^N.$ 
	\end{cor}
	
	\begin{proof}
	This follows from the $k=2$ case of Theorem \ref{thm:TopBound}, which implies that
	
		\begin{equation*}
			\bigg{|} 1 - \frac{I_N}{e^{N^2F_N^{(0)} + F_N^{(1)} + N^{-2}F_N^{(2)}}} \bigg{|} <1
			\quad\text{ and }\quad
			\bigg{|} 1 - \frac{J_N}{e^{N^2G_N^{(0)} + G_N^{(1)} + N^{-2}G_N^{(2)}}} \bigg{|} <1
		\end{equation*}
				
	\noindent
	on $\overline{\D}_{\gamma}^N$ for $N$ sufficiently large. These inequalities in turn imply the 
	non-vanishing of 
	
		\begin{equation*}
			\frac{I_N}{e^{N^2F_N^{(0)} + F_N^{(1)} + N^{-2}F_N^{(2)}}}
			\quad\text{ and }\quad
			\frac{J_N}{e^{N^2G_N^{(0)} + G_N^{(1)} + N^{-2}G_N^{(2)}}}
		\end{equation*}

	\noindent
	on $\overline{\D}_\gamma^N,$ from which we conclude the nonvanishing of $I_N$ and $J_N$ on this polydisc.
	\end{proof}
	
	\subsection{Analytic error functions}
	Set $\xi = \min(\gamma,\delta),$ where $\gamma$ is the positive constant in 
	Theorem \ref{thm:TopBound} and $\delta$ is the positive constant in Theorem \ref{thm:AnalyticCandidates}.
	We may define an array of analytic functions on $\overline{\D}_\xi^N$ by
	
	\begin{align*}
		\Delta_N^{(0)} &= N^{-2}\log I_N - F_N^{(0)} \\
		\Delta_N^{(k)} &= N^2\Delta_N^{(k-1)} - F_N^{(k)}, \quad k \in \mathbb{N}.
	\end{align*}

	\noindent
	Explicitly, we have

	\begin{equation*}
		\Delta_N^{(k)} = N^{2k-2} \bigg{(} \log I_N - \sum_{g=0}^k N^{2-2g} F_N^{(g)} \bigg{)}, \quad k \in \mathbb{N}
		\cup \{0\}.
	\end{equation*}
	
	\noindent
	We could also have defined $\Delta_N^{(k)}$ using $J_N$ in place of $I_N$, and $G_N^{(g)}$ in place
	of $F_N^{(g)}$, and in what follows $\Delta_N^{(k)}$ may equally well be replaced with this function instead.
	Our main result, Theorem \ref{thm:Main}, is an immediate consequence of the following 
	convergence theorem, the proof of which occupies the remainder of the paper.

	\begin{thm}
	\label{thm:MainReformulated}
		For any  $0 < \varepsilon < \xi$ 
		we have $\lim_{N \to \infty} \| \Delta_N^{(k)} \|_\varepsilon = 0$
		for each $k \in \mathbb{N}$.
	\end{thm}
	
	\subsection{Reduction to uniform boundedness}
	By virtue of its definition, the function $\Delta_N^{(k)}$ admits the string expansion
	
		\begin{equation*}
			\Delta_N^{(k)} = \sum_{d=1}^\infty \frac{z^d}{d!} \sum_{\alpha,\beta \vdash d}
			\frac{p_\alpha(a_1,\dots,a_N)}{N^{\ell(\alpha)}} \frac{p_\beta(b_1,\dots,b_N)}{N^{\ell(\beta)}}
			\Delta_N^{(k)}(\alpha,\beta),
		\end{equation*}
		
	\noindent
	the coefficients of which are given by
	
		\begin{equation*}
			\Delta_N^{(k)}(\alpha,\beta) = N^{2k-2}\left( L_N(\alpha,\beta) - 
			\sum_{g=0}^k \frac{\mon_g(\alpha,\beta)}{N^{2g}} \right),
		\end{equation*}
		
	\noindent
	where 
	
		\begin{equation*}
			\log I_N = \sum_{d=1}^\infty \frac{z^d}{d!} \sum_{\alpha,\beta \vdash d} 
			\frac{p_\alpha(a_1,\dots,a_N)}{N^{\ell(\alpha)}} \frac{p_\beta(b_1,\dots,b_N)}{N^{\ell(\beta)}}
			L_N(\alpha,\beta),
		\end{equation*}
		
	\noindent
	and both series converge absolutely on $\overline{\D}_\gamma^N.$ From Section \ref{sec:Stable},
	we know that, for each $k \in \mathbb{N} \cup \{0\},$ each fixed string coefficient of $\Delta_N^{(k)}$
	converges to zero as $N \to \infty,$
	
		\begin{equation*}
			\lim_{N \to \infty} \Delta_N^{(k)}(\alpha,\beta) = 0.
		\end{equation*}
		
	\noindent
	In fact, asymptotic vanishing of the string coefficients of $\Delta_N^{(k)}$ implies 
	uniform asymptotic vanishing of string series provided we have uniform boundedness.

	Let $m \in \mathbb{N}$ be an arbitrary positive integer. By a ``normalized string series'' on 
	
		\begin{equation*}
			\overline{\D}_\xi \times \overline{\D}_1^{mN},
		\end{equation*} 
		
	\noindent
	we mean a power series of the form

	\begin{equation*}
		\Delta_N = \sum_{d=1}^\infty \frac{z^d}{d!} \sum_{\alpha^1,\dots,\alpha^m \vdash d}
		\prod_{i=1}^m \frac{p_{\alpha^i}(a_{i1},\dots,a_{iN})}{N^{\ell(\alpha^i)}}\Delta_N(\alpha^1,\dots,\alpha^m) 
	\end{equation*}
	
	\noindent
	which converges absolutely on $\overline{\D}_\xi^N$. In order to 
	prove Theorem \ref{thm:MainReformulated}, we will use the fact that, in the presence of uniform boundedness, 
	uniform convergence of $\Delta_N$ on any closed proper subset of $\overline{\D}_\xi^N$ follows from the
	convergence of each of its string coefficients $\Delta_N(\alpha^1,\dots,\alpha^m)$.

	\begin{lem}
		\label{lem:ConvergenceLemma}
		If  $\sup_{N \in \mathbb{N}} \|\Delta_N\|_\xi < \infty$ and 
			
			\begin{equation*}
				\lim_{N \to \infty} \Delta_N(\alpha^1,\dots,\alpha^m)=0
			\end{equation*}
			
		\noindent
		for any $d \in \mathbb{N}$ and $\alpha^1,\dots,\alpha^m \vdash d$, then
		
			\begin{equation*}
				\lim_{N \to \infty} \|\Delta_N\|_\varepsilon = 0
			\end{equation*}
			
		\noindent
		for any $0 < \varepsilon < \xi$.
	\end{lem}
	
	\begin{proof}
	Fix $\varepsilon < \xi.$ Let $\kappa > 0$ be given. For any $n,N \in \mathbb{N}$, 			
			we have
			
				\begin{align*}
					\left\| \Delta_N \right\|_\varepsilon
					\leq& \sum_{d=1}^n \frac{\varepsilon^d}{d!}\sum_{\alpha^1,\dots,\alpha^m \vdash d}  
					\left| \Delta_N(\alpha^1,\dots,\alpha^m)  \right|
					 \\ 
					+& 
					\sum_{d=n+1}^\infty \frac{\varepsilon^d}{d!} \left| \sum_{\alpha^1,\dots,\alpha^m \vdash d} 
					 \Delta_N(\alpha^1,\dots,\alpha^m) 
					\prod_{i=1}^m \frac{p_{\alpha^i}(a_1,\dots,a_N)}{N^{\ell(\alpha^i)}} \right|,
				\end{align*}
				
			\noindent
			by the triangle inequality. By Cauchy's estimate, 
			
				\begin{equation*}
					\frac{1}{d!} \left| \sum_{\alpha^1,\dots,\alpha^m \vdash d} 
					 \Delta_N(\alpha^1,\dots,\alpha^m) 
					\prod_{i=1}^m \frac{p_{\alpha^i}(a_1,\dots,a_N)}{N^{\ell(\alpha^i)}} \right|
					\leq \frac{\|\Delta_N\|_\xi}{\xi^d}.
				\end{equation*}
				
			\noindent
			Thus
			
				\begin{equation*}
					\left\| \Delta_N \right\|_{\varepsilon} \\
					\leq \ \sum_{d=1}^n \frac{\varepsilon^d}{d!}\sum_{\alpha^1,\dots,\alpha^m \vdash d}  
					\left| \Delta_N(\alpha^1,\dots,\alpha^m)  \right|
					 + \left( \frac{\varepsilon}{\xi} \right)^{n+1} K,
				\end{equation*}

			\noindent
			where 
			
				\begin{equation*}
					K = \frac{\sup_{N \in \mathbb{N}} \|\Delta_N\|_\xi}{1-\frac{\varepsilon}{\xi}}
				\end{equation*}
				
			\noindent
			is a constant. Since 
			
				\begin{equation*}
					\lim_{n \to \infty} \left( \frac{\varepsilon}{\xi} \right)^{n+1} = 0,
				\end{equation*}
				
			\noindent
			we can choose $n_0$ sufficiently large so that 
			
				\begin{equation*}
					\left( \frac{\varepsilon}{\xi} \right)^{n_0+1} K < \frac{\kappa}{2}.
				\end{equation*}
				
			\noindent
			Then, since 
				
				\begin{equation*}
					\lim_{N \to \infty} |\Delta_N(\alpha^1,\dots,\alpha^m)| = 0
				\end{equation*}
				
			\noindent
			for each $d \in \mathbb{N}$ and all $\alpha^1,\dots,\alpha^m \vdash d$, we can choose $N_0$ 
			sufficiently large so that $N \geq N_0$ implies 
			
				\begin{equation*}
					 \sum_{d=1}^{n_0} \frac{\varepsilon^d}{d!}\sum_{\alpha^1,\dots,\alpha^m \vdash d}  
					 \left| \Delta_N(\alpha^1,\dots,\alpha^m)  \right|
					  < \frac{\kappa}{2}.
				\end{equation*}
				
			\noindent
			We conclude that $N \geq N_0$ implies
			
				\begin{equation*}
					\left\| \Delta_N \right\|_{\varepsilon}  < \kappa,
				\end{equation*}
				
			\noindent
			as required.

	\end{proof}		
	
	\subsection{Proof of uniform boundedness}
	In view of Lemma \ref{lem:ConvergenceLemma}, the following 
	result completes the proof of Theorem \ref{thm:MainReformulated} 
	and hence also of Theorem \ref{thm:Main}, which proves Conjecture \ref{conj:Main}.
			
	\begin{thm}
		\label{thm:UniformBoundedness}
		For any $\varepsilon < \xi,$ we have
		
			\begin{equation*}
				\sup_{N \in \mathbb{N}} \|\Delta_N^{(k)}\|_\varepsilon < \infty
			\end{equation*}
			
		\noindent
		for each $k \in \mathbb{N} \cup \{0\}.$
	\end{thm}
		
	\begin{proof}
	Let $(\varepsilon_k)_{k=0}^\infty$ be a strictly decreasing sequence in the 
	interval $(\varepsilon,\xi).$ Then, for each $N \in \mathbb{N}$, we have a 
	corresponding sequence of nested closed polydiscs,
	
		\begin{equation*}
			\overline{\D}_\xi^N \supset \overline{\D}_{\varepsilon_0}^N \supset \overline{\D}_{\varepsilon_1}^N
			\supset \dots \supset  \overline{\D}_{\varepsilon}^N.
		\end{equation*}	
	
		Let $k \in \mathbb{N} \cup \{0\}$ be fixed. Observe that 
		
			\begin{equation*}
				\frac{I_N}{e^{\sum_{g=0}^k N^{2-2g}F_N^{(g)}}} = e^{N^{2-2k}\Delta_N^{(k)}}.
			\end{equation*}		
			
		\noindent
		Thus, by Theorem \ref{thm:TopBound}, we have 
		
			\begin{equation*}
				\bigg{|} 1 - e^{N^{2-2k}\Delta_N^{(k)}} \bigg{|} \leq c_k N^{2-2k}
			\end{equation*}
			
		\noindent
		on $\overline{\D}_{\varepsilon_k}^N$ for $N$ sufficiently large, where $c_k$ is a positive
		constant depending only on $k.$ This in turn implies

			\begin{equation*}
				\left| e^{N^{2-2k}\Delta_N^{(k)}} \right| \leq 1 + c_kN^{2-2k}
			\end{equation*}
			
		\noindent
		on $\overline{\D}_{\varepsilon_k}^N$ for $N$ sufficiently large.
		Thus
		
			\begin{equation*}
				N^{2-2k} \operatorname{Re} \Delta_N^{(k)} \leq \log \left(1 + c_kN^{2-2k} \right).
			\end{equation*}
			
		\noindent
		on $\overline{\D}_{\varepsilon_k}^N$ for $N$ sufficiently large. If $k \neq 1,$ this yiels
		
			\begin{equation*}
				N^{2-2k} \operatorname{Re} \Delta_N^{(k)} \leq \log \left(1 + c_kN^{2-2k} \right) \leq c_kN^{2-2k}
			\end{equation*}
			
		\noindent
		on $\overline{\D}_{\varepsilon_k}^N$ for $N$ sufficiently large, 
		so that 
		
			\begin{equation*}
				\operatorname{Re} \Delta_N^{(k)} \leq c_k
			\end{equation*}
			
		\noindent
		on $\overline{\D}_{\varepsilon_k}^N$ for $N$ sufficiently large.
		If $k=1$ we obtain instead
		
			\begin{equation*}
				\operatorname{Re} \Delta_N^{(k)} \leq \log \left(1 + c_k \right)
			\end{equation*}
			
		\noindent
		on $\overline{\D}_{\varepsilon_k}^N$ for $N$ sufficiently large.
		Thus in all cases we have a uniform bound on the real part of $\Delta_N^{(k)}$, i.e.
		a bound of the form
		
			\begin{equation*}
				\operatorname{Re} \Delta_N^{(k)} \leq \tilde{c}_k
			\end{equation*}

		\noindent
		for some positive constant $\tilde{c}_k.$ In order to leverage this into a bound
		on the modulus, we apply the Borel-Carath\'eodory inequality (see e.g. \cite{Titchmarsh}), which 
		bounds the sup norm of an analytic function on a closed disc in terms of the supremum of its real
		part on a larger closed disc. Applying the Borel-Carath\'eodory inequality, we obtain the bound
				
			\begin{equation*}
				\| \Delta_N^{(k)} \|_{\varepsilon_{k+1}} \leq 
				\frac{2\varepsilon_{k+1}}{\varepsilon_k - \varepsilon_{k+1}}\ \sup_{\overline{\D}_{\varepsilon_k}}\
				\operatorname{Re} \Delta_N^{(k)} \leq 
				\frac{2\varepsilon_{k+1}}{\varepsilon_k-\varepsilon_{k+1}} \tilde{c}_k.
			\end{equation*}
			
		\noindent
		Since $\varepsilon < \varepsilon_k,$ we have 
		
			\begin{equation*}
				\|\Delta_N^{(k)} \|_{\varepsilon} \leq\|\Delta_N^{(k)} \|_{\varepsilon_{k+1}} \leq
				\frac{2\varepsilon_{k+1}}{\varepsilon_k-\varepsilon_{k+1}} \tilde{c}_k,
			\end{equation*}
			
		\noindent
		as required.
	\end{proof} 	
	
		 
\bibliographystyle{amsplain}

\end{document}